\documentclass[10pt,reqno]{amsart}
\usepackage{kotex}
\usepackage[T1]{fontenc}
\usepackage[utf8]{inputenc}
\usepackage{lmodern}
\usepackage{amsmath,amssymb,amsthm,mathtools,mathrsfs}
\usepackage{microtype}
\usepackage{xcolor}
\usepackage{aliascnt}
\usepackage[pdfusetitle,colorlinks=true,linkcolor=blue!55!black,
  citecolor=green!40!black,urlcolor=blue!60!black]{hyperref}
\usepackage[nameinlink,capitalize,noabbrev]{cleveref}

\allowdisplaybreaks
\numberwithin{equation}{section}

\newtheorem{theorem}{Theorem}[section]
\newaliascnt{proposition}{theorem}
\newtheorem{proposition}[proposition]{Proposition}
\aliascntresetthe{proposition}
\newaliascnt{lemma}{theorem}
\newtheorem{lemma}[lemma]{Lemma}
\aliascntresetthe{lemma}
\newaliascnt{corollary}{theorem}
\newtheorem{corollary}[corollary]{Corollary}
\aliascntresetthe{corollary}
\theoremstyle{definition}
\newaliascnt{definition}{theorem}
\newtheorem{definition}[definition]{Definition}
\aliascntresetthe{definition}
\newaliascnt{example}{theorem}

\aliascntresetthe{example}
\theoremstyle{remark}
\newaliascnt{remark}{theorem}
\newtheorem{remark}[remark]{Remark}
\aliascntresetthe{remark}

\crefname{theorem}{Theorem}{Theorems}
\crefname{proposition}{Proposition}{Propositions}
\crefname{lemma}{Lemma}{Lemmas}
\crefname{corollary}{Corollary}{Corollaries}
\crefname{definition}{Definition}{Definitions}
\crefname{example}{Example}{Examples}
\crefname{remark}{Remark}{Remarks}
\Crefname{theorem}{Theorem}{Theorems}
\Crefname{proposition}{Proposition}{Propositions}
\Crefname{lemma}{Lemma}{Lemmas}
\Crefname{corollary}{Corollary}{Corollaries}
\Crefname{definition}{Definition}{Definitions}
\Crefname{example}{Example}{Examples}
\Crefname{remark}{Remark}{Remarks}

\newcommand{\Hh}{\mathbb H}
\newcommand{\Cc}{\mathbb C}
\newcommand{\Rr}{\mathbb R}
\newcommand{\Zz}{\mathbb Z}
\newcommand{\Gam}{\Gamma}
\newcommand{\FP}{\operatorname*{FP}}
\newcommand{\Id}{\operatorname{Id}}
\newcommand{\cR}{\mathcal R}
\newcommand{\cL}{\mathcal L}
\newcommand{\Resol}{\mathscr R}
\newcommand{\cU}{\mathbf U}
\newcommand{\cK}{\mathbf K}
\newcommand{\Kdir}{\mathcal K}
\newcommand{\Dkt}{\mathcal D_{k,t}}
\newcommand{\lam}{\lambda}
\newcommand{\nK}{\mathfrak n}
\newcommand{\Poch}[2]{(#1)_{#2}}

\title{Weight-Shifting Resolvent Kernels on the Modular Surface}

\author{SEUNGJU LEE}
\date{\today}
\address{Department of Mathematics, Ajou University, 206 Worldcup-ro,
Yeongtong-gu, Suwon 16499, Republic of Korea}
\email{chesohater@ajou.ac.kr}


\begin{document}

\begin{abstract}
The local hypergeometric resolvent kernels and their mixed \(K\)-type
shifts are due to Fay.  On the modular surface, we determine the
normalization relating these
kernels to the Maass-shifted automorphic resolvent family and equip that
family with two compatible extensions: the hyperbolic circular finite part at
the diagonal and a parameter-dependent Mellin subtraction matched to the
cusp zero mode.  On moderate-growth eigenfunctions the resulting meromorphic transform
satisfies an explicit spectral-reduction formula.  Suppose that the
multiplier is trivial at the cusp and that such an eigenfunction has
spectral parameter \(s_0\ne1/2\).  If \(s_0\) is a simple zero of the
diagonal source coefficient, the continued kernel and Eisenstein series
are regular at \(s_0\), and the specialized integral converges, then
specializing the kernel before integration need not agree with continuing
the parameter-dependent finite-part family.  Their difference is the
relevant constant-term coefficient of the input multiplied by the first
spectral-parameter derivative of the kernel's constant Fourier coefficient.
For the trivial multiplier, even \(k>2\), raising shifts \(t=k+2r\), and
asymptotic normalization \((k-1)/(4\pi)\), the holomorphic point is
\(s=k/2\): specialized integration depends only on the cuspidal part of a
modular form, whereas continued finite-part integration acts on the full
form.  In the case \(k=12\) and \(t=14\), the specialized kernel annihilates
\(y^6E_{12}^{\mathrm{hol}}\), while the continued value is
\(-\cR_{12}(y^6E_{12}^{\mathrm{hol}})/12\), where
\(\cR_{12}=2iy\partial_z+6\) is the Maass raising operator from weight
\(12\) to weight \(14\).  We also record the Maass shift of the spectral projection and show that
the value and first spectral derivative together recover the Maass shift
of the reduced resolvent.
\end{abstract}

\maketitle

\section{Introduction and statement of results}
\label{sec:introduction}

A meromorphic automorphic kernel on a noncompact surface presents two
distinct limiting problems.  Its diagonal singularity must be extended
as a distribution, while its constant Fourier mode may require a cusp
subtraction.  A further choice arises at a special spectral parameter:
one may first specialize the kernel and then integrate, or first form a
meromorphic family of finite-part integrals and then continue to the
special parameter.  These operations need not commute, even when the
already specialized integral converges.  The most concrete instance in
this paper is
\begin{equation}
 \cU_{6,\mathrm{sp}}^{14,12}
   \bigl(y^6E_{12}^{\mathrm{hol}}\bigr)=0,
 \qquad
 \operatorname{AC}_{s=6}\cU_{s,\mathrm{fp}}^{14,12}
   \bigl(y^6E_{12}^{\mathrm{hol}}\bigr)
 =-\frac1{12}\cR_{12}\bigl(y^6E_{12}^{\mathrm{hol}}\bigr)\ne0.
 \label{eq:intro-e12}
\end{equation}
Here the subscript \(\mathrm{sp}\) denotes ordinary integration after
specializing the kernel, whereas \(\operatorname{AC}\) denotes the
removable value of the parameter-dependent finite-part family.  Thus the
second expression is not another regularized value assigned to the fixed
integrable function in the first expression.

We work first with a unitary multiplier of real weight \(k\) on
\(\Gam=\mathrm{SL}_2(\Zz)\).  If \(t-k\in2\Zz\), let \(\Dkt\) be the
corresponding iterated Maass raising or lowering operator and let
\(\Resol_k(s)=(\Delta_k-s(1-s))^{-1}\).  The central family is
\begin{equation}
 \cU_s^{t,k}=c_{t,k}(s)\Dkt\Resol_k(s),\qquad
 c_{t,k}(s)=4\pi\nK(-1)^r
 \frac{\Gamma(2s)}
 {\Gamma(s-\ell_-/2)\Gamma(s+\ell_+/2)},
 \label{eq:intro-family}
\end{equation}
where \(\nK\in\Cc^\times\) is a fixed nonzero asymptotic normalization,
\(r=|t-k|/2\), \(\ell_-=\min\{k,t\}\), and
\(\ell_+=\max\{k,t\}\).  The scalar \(c_{t,k}\) is the diagonal source
coefficient of the explicit mixed kernel.  The product in \eqref{eq:intro-family} is first defined for
\(\Re s>1/2\) when \(s(1-s)\) lies in the resolvent set of \(\Delta_k\),
and then continued meromorphically.

The local ingredients in \eqref{eq:intro-family} are known.  For a generic
spectral parameter, Fay's Theorem~1.1 and Equation~(13) give a fundamental
pair for the radial equation in each angular Fourier mode.  For the singular
\(n=0\) member used as the resolvent seed, the prescribed subdominant
asymptotic selects a one-dimensional line, and the logarithmic diagonal
coefficient fixes its scalar once the unit source is specified
\cite[Theorem~1.1 and Equations~(13), (17); see also the proof of
Theorem~1.4]{Fay77}.  Fay's Theorem~2.1 and
Equations~(45)--(46) then transport this normalized solution to the mixed
\(K\)-types \cite[Theorem~2.1 and Equations~(45)--(46)]{Fay77}.
Thus covariance and the eigen-equations alone are not the uniqueness statement.  Kara et al.
give a convenient real-weight automorphic resolvent realization
\cite[Equation~(16)]{KaraEtAl21}, and related differential identities
appear in \cite[Theorem~2.6]{BringmannJenkinsKane19}.  Fay's Equation~(39) and Corollary~3.2 exhibit the residue and
reproducing-kernel mechanism at the holomorphic point~\cite[Equation~(39) and Corollary~3.2]{Fay77};
compare the reproducing-kernel identity in
\cite[Definition~7 and Theorem~4(3)]{Steiner16}.
Accordingly, we
do not claim a new hypergeometric kernel, a new mixed-shift identity, or
a new Bergman projector.

Our task is to define the action of that known kernel consistently at
both the diagonal and the cusp.  For higher shifts \(r\geq2\) the
mixed kernel is generally not locally integrable; for \(r=0,1\) the same
prescription agrees with ordinary local integration.  We use the circular
finite part in a lifted hyperbolic disk coordinate.  For \(z,w\in X\)
sufficiently close, choose a lift \(\widetilde z\in\Hh\) of \(z\) and the
compatible local lift \(\widetilde w\in\Hh\) of \(w\), and set
\[X_{\widetilde z}(\widetilde w)=\frac{\widetilde w-\widetilde z}{\widetilde w-\overline{\widetilde z}}.\]
Changing the compatible lifts by
\(\gamma=\left(\begin{smallmatrix}a&b\\c&d\end{smallmatrix}\right)\in\Gam\)
multiplies this coordinate by
\((c\overline{\widetilde z}+d)/(c\widetilde z+d)\), which has modulus one;
hence the circular cutoff is independent of the lifts.  The resulting
extension still depends on the cutoff geometry: changing the shape of the
small deleted disk changes the extension by a distribution supported on
the diagonal.  A local calculation shows that iterated raising or
iterated lowering of the logarithmic Green singularity produces this
circular extension and no additional diagonal-supported term in the
kernel distribution.  When the cusp character is trivial, the continued
kernel has constant term
\begin{equation}
 \beta_s(z)v^{1-s},\qquad
 \beta_s(z)=\frac{c_{t,k}(s)}{2s-1}\Dkt E_k(z,s).
 \label{eq:intro-beta}
\end{equation}
We subtract the corresponding Mellin antiderivative of its product with the two
constant exponents of the input. This follows the truncation--Mellin framework
of Zagier
\cite[Equations~(15), (27), and~(29)--(31)]{Zagier81}.  In a related setting, Wu shows that taking the residue at $s=\tfrac12$
and letting the cusp-truncation parameter $T$ tend to infinity need not
commute \cite[Section~1.3 and Theorem~2.12(2), (4)]{Wu19};
here the objective is the resulting defect for the Fay family.

With both extensions normalized, the main calculation reduces to the
following identity.  If
\(\Delta_k\phi=\lambda_\phi\phi\) and \(\phi\) has moderate growth, then
\begin{equation}
  {\;
 \cU_{s,\mathrm{fp}}^{t,k}\phi
 =\frac{c_{t,k}(s)}
 {\lambda_\phi-s(1-s)}\Dkt\phi .
 \;}
 \label{eq:intro-spectral-reduction}
\end{equation}
This is \cref{thm:regularized-reduction}.  It is proved by Green's
identity after matching the cusp counterterm with the boundary Wronskian,
without
differentiating a moving diagonal cutoff.

When the cusp character is trivial, the main general consequence is
\cref{thm:boundary-defect}.  Suppose
\(s_0\ne1/2\), \(c_{t,k}\) has a simple zero at \(s_0\), the continued
kernel and \(E_k(\cdot,s)\) are regular there, and a moderate-growth
eigenfunction has constant term
\[
 A_\phi y^{s_0}+B_\phi y^{1-s_0}.
\]
If its specialized kernel integral converges, then
\begin{equation}
  {\;
 \operatorname{AC}_{s=s_0}\cU_{s,\mathrm{fp}}^{t,k}\phi
 -\cU_{s_0,\mathrm{sp}}^{t,k}\phi
 =A_\phi\frac{c_{t,k}'(s_0)}{2s_0-1}
   \Dkt E_k(\cdot,s_0).
 \;}
 \label{eq:intro-boundary-defect}
\end{equation}
Thus the discrepancy is neither arbitrary nor a competing convention for
a fixed integral: it is \(A_\phi\) times the first \(s\)-derivative of the
coefficient in \eqref{eq:intro-beta}.

The holomorphic specialization makes
\eqref{eq:intro-boundary-defect} explicit.  Let \(k>2\) be even,
\(t=k+2r\), \(\nK=(k-1)/(4\pi)\), and
\[
 F^\circ=F-a_0(F)E_k^{\mathrm{hol}},\qquad
 \Phi_F=y^{k/2}F,\qquad
 C_{k,r}=\frac{(-1)^r}{(k)_r}.
\]
For every \(F\in M_k(\Gam)\), \cref{thm:holomorphic-defect} gives
\begin{equation}
  {\;
 \cU_{k/2,\mathrm{sp}}^{k+2r,k}\Phi_F
 =C_{k,r}\cR_k^{(r)}\Phi_{F^\circ},\qquad
 \operatorname{AC}_{s=k/2}\cU_{s,\mathrm{fp}}^{k+2r,k}\Phi_F
 =C_{k,r}\cR_k^{(r)}\Phi_F .
 \;}
 \label{eq:intro-holomorphic-defect}
\end{equation}
Their difference is the raised Eisenstein term proportional to \(a_0(F)\).
\Cref{cor:weight-twelve-defect} is the case displayed in
\eqref{eq:intro-e12}; it proves that the boundary contribution can be
nonzero while the ordinary specialized integral vanishes.

For \(\Re s_0>1/2\), a second consequence concerns a simple zero of
\(c_{t,k}\) matched to an isolated spectral point.  The reduced-resolvent
expansion \cite[Chapter~III, Section~6.5, especially Equation~(6.32)]{Kato95} gives
\[\cU_{s_0}^{t,k}=\frac{c_{t,k}'(s_0)}{2s_0-1}\Dkt P_0,\]
and the value and first derivative of the explicit kernel family together
recover \(\Dkt S_0\), the Maass shift of the reduced resolvent; see
\cref{thm:selector}.  This follows from the Laurent expansion of the
resolvent at an isolated eigenvalue.  Its role here is to identify which
spectral operators are represented by the
normalized Fay kernels and their first parameter derivatives.

The paper is confined to the modular surface and its single cusp.  The
finite-part normalization is fixed by the hyperbolic circular cutoff, the
zero-mode Mellin antiderivative, and agreement with the kernel integral in
\eqref{eq:direct-agreement} where it converges.  No moving-cutoff
regularity theorem is asserted.  \Cref{sec:fay-kernels} establishes the
Fay normalization dictionary; \cref{sec:direct-kernels} constructs the
diagonal extension; \cref{sec:continuation} computes the cusp coefficient;
\cref{sec:regularized-action} proves
\eqref{eq:intro-spectral-reduction}; and
\cref{sec:spectral-selectors} proves the selector and specialization
defect formulas.

\section{Fay kernels and normalized Maass shifts}
\label{sec:fay-kernels}

\subsection{Conventions}

Let \(\Gam=\mathrm{SL}_2(\Zz)\), let \(\overline\Gam=\Gam/\{\pm I\}\), put \(X=\overline\Gam\backslash\Hh\), and write \(z=x+iy\in\Hh\). The hyperbolic measure and the weight-\(k\) Laplacian are
\begin{equation} d\mu(z)=\frac{dx\,dy}{y^2},\qquad \Delta_k=-y^2(\partial_x^2+\partial_y^2)+iky\partial_x. \label{eq:laplacian}\end{equation}

For \(\gamma=\left(\begin{smallmatrix}a&b\\c&d\end{smallmatrix}\right)\) put
\[j_k(\gamma,z)=\left(\frac{cz+d}{|cz+d|}\right)^k,\]
with the principal argument. We fix a unitary multiplier system \(\nu\) of real weight \(k\). If \(t-k\in2\Zz\), the same multiplier may be used at weight \(t\). Thus our automorphic functions satisfy
\[f(\gamma z)=\nu(\gamma)j_k(\gamma,z)f(z).\]
The usual cocycle qualifications for nonintegral weights are understood; see \cite[Sections~2.2--2.3]{KaraEtAl21}, with their weight index equal to \(k/2\) in our convention. They play no role in the local calculations below.

Let \(L_k^2(X,\nu)\) denote the Hilbert space of measurable functions with this transformation law and finite norm
\[\|f\|_2^2=\int_{\overline\Gam\backslash\Hh}|f(z)|^2\,d\mu(z).\]
We write \(C_{\mathrm{c},k}^\infty(X,\nu)\) for the space of smooth functions satisfying the above weight-\(k\) transformation law and having compact support in \(X\), and \(C_{\mathrm{loc},k}^\infty(X,\nu)\) for the space of all smooth functions satisfying the same transformation law, endowed with the topology of \(C^\infty\)-convergence on compact subsets of \(X\). A kernel of weight \(t\) in \(z\) and dual weight \(-k\) in \(w\) has multiplier \(\nu\) in the first variable and \(\bar\nu\) in the second; consequently \(K(z,w)f(w)d\mu(w)\) descends to the quotient.

On smooth weight-\(\ell\) functions define the Maass operators
\begin{equation} \cR_\ell=2iy\partial_z+\frac{\ell}{2},\qquad \cL_\ell=-2iy\partial_{\bar z}-\frac{\ell}{2}. \label{eq:maass-operators}\end{equation}
They raise and lower the weight by two and satisfy
\begin{equation} \Delta_{\ell+2}\cR_\ell=\cR_\ell\Delta_\ell,\qquad \Delta_{\ell-2}\cL_\ell=\cL_\ell\Delta_\ell. \label{eq:maass-intertwining}\end{equation}
For \(r\in\Zz_{\geq0}\) set
\[\cR_k^{(r)}=\cR_{k+2r-2}\cdots\cR_k,\qquad \cL_k^{(r)}=\cL_{k-2r+2}\cdots\cL_k,\]
where the empty products are the identity. If
\[m=\frac{t-k}{2}\in\Zz,\qquad r=|m|,\]
we use the single notation
\begin{equation} \Dkt=\begin{cases}\cR_k^{(r)},&t>k,\\ \Id,&t=k,\\ \cL_k^{(r)},&t<k.\end{cases} \label{eq:shift-operator}\end{equation}
Then \(\Delta_t\Dkt=\Dkt\Delta_k\). We reserve \(\Resol_k(s)\), rather than \(R_k\), for the resolvent, thereby avoiding an ambiguity between a resolvent and a raising operator.

Put
\begin{equation} \lam(s)=s(1-s),\qquad \ell_-=\min\{k,t\},\qquad \ell_+=\max\{k,t\}. \label{eq:spectral-parameter}\end{equation}
The standard self-adjoint realization of \(\Delta_k\) on \(L_k^2(X,\nu)\) is denoted by the same symbol. Whenever \(\lam(s)\) lies in its resolvent set, put
\begin{equation} \Resol_k(s)=(\Delta_k-\lam(s))^{-1}. \label{eq:resolvent-definition}\end{equation}
We work initially in the half-plane \(\Re s>1/2\). For \(s\) in this half-plane with \(\lam(s)\) in the resolvent set, \eqref{eq:resolvent-definition} is the \(L^2\)-resolvent.

\subsection{The explicit seed and Fay's normalization}

For distinct \(z=x+iy,w=u+iv\in\Hh\), let
\[U(z,w)=\frac{|z-w|^2}{4yv},\qquad \sigma(z,w)=1+U(z,w).\]
Fix a nonzero asymptotic normalization \(\nK\in\Cc^\times\), and define
\begin{align}
 P_{t,k}(z,w)&=(yv)^{t/2}\left(\frac{2i}{z-\bar w}\right)^t\left(\frac{\bar w-z}{w-z}\right)^m,\label{eq:covariant-factor}\\
 K_s^{t,k}(z,w)&=\nK P_{t,k}(z,w)U(z,w)^{t/2-s}{}_2F_1\!\left(s-\frac t2,s-\frac k2;2s;-\frac1{U(z,w)}\right).\label{eq:explicit-seed}
\end{align}
The last power in \eqref{eq:covariant-factor} is integral. Since \(2i/(z-\bar w)\) lies in the right half-plane, its real power is unambiguous with the principal logarithm. The coefficient of the prescribed asymptotic at hyperbolic infinity is \(\nK\); the corresponding diagonal source coefficient is computed below.

The equal-weight seed is Fay's kernel with a different automorphy-factor convention and normalization.

\begin{proposition}[Fay--Kara normalization dictionary]
\label{prop:fay-dictionary}
Let
\[H_{k/2}(z,w)=\left(-\frac{|z-\bar w|^2}{(z-\bar w)^2}\right)^{k/2}.\]
The power is defined using the principal logarithm.  For \(t=k\), Pfaff's transformation gives
\begin{equation} K_s^{k,k}(z,w)=\nK H_{k/2}(z,w)\sigma(z,w)^{-s}{}_2F_1\!\left(s-\frac k2,s+\frac k2;2s;\sigma(z,w)^{-1}\right). \label{eq:equal-seed-pfaff}\end{equation}
If \(G_{s,k}^{\mathrm{Fay}}(z,w)\) denotes the unit-source seed in our weight convention, then
\begin{equation} K_s^{k,k}=c_{k,k}(s)G_{s,k}^{\mathrm{Fay}},\qquad c_{k,k}(s)=4\pi\nK\frac{\Gamma(2s)}{\Gamma(s-k/2)\Gamma(s+k/2)}. \label{eq:fay-identification}\end{equation}
Kara et al.'s weight index is \(k/2\) when the unitary weight used here is \(k\).
\end{proposition}

\begin{proof}
Since \(U=\sigma-1\), Pfaff's identity \cite[Equation~15.8.1]{DLMF15} yields
\[U^{-a}{}_2F_1(a,a;2s;-U^{-1})=\sigma^{-a}{}_2F_1(a,2s-a;2s;\sigma^{-1}),\qquad a=s-\frac k2.\]
Equation \eqref{eq:covariant-factor} gives \(P_{k,k}=\sigma^{-k/2}H_{k/2}\), proving \eqref{eq:equal-seed-pfaff}. Kara et al.'s point-pair seed, displayed after their automorphic kernel formula \cite[Equation~(16) and the displayed definition following it]{KaraEtAl21}, is
\[H_{k/2}\sigma^{-s}\frac{\Gamma(s-k/2)\Gamma(s+k/2)}{4\pi\Gamma(2s)}{}_2F_1(s+k/2,s-k/2;2s;\sigma^{-1}).\]
Comparison gives \eqref{eq:fay-identification}. Fay defines the radial hypergeometric solutions in \cite[Equation~(13)]{Fay77}. In the proof of Theorem~1.4, he constructs the fundamental solution on \(\Hh\) from the singular \(n=0\) radial function and identifies it as the resolvent, while Equation~(29) gives the corresponding spectral inversion formula \cite[Theorem~1.4 and Equation~(29)]{Fay77}. After translating the weight and automorphy-factor conventions, that fundamental solution agrees with the seed above. With Fay's operator convention \cite[p.~144]{Fay77}, the operator at index \(k/2\) is \(D_{k/2}^{\mathrm{Fay}}=-\Delta_k\), so
\begin{equation} (\Delta_k-s(1-s))^{-1}=-\bigl(D_{k/2}^{\mathrm{Fay}}-s(s-1)\bigr)^{-1}. \label{eq:fay-laplacian-sign}\end{equation}
This accounts for the sign change between Fay's resolvent and the unit-source convention in \eqref{eq:resolvent-definition}.
\end{proof}

\begin{proposition}[Normalized mixed shifts]
\label{prop:normalized-shifts}
For \(r\in\Zz_{\geq0}\), the following identities hold where the Pochhammer denominators are nonzero and thereafter as meromorphic identities in \(s\):
\begin{align}
 K_s^{k+2r,k}&=\frac{(-1)^r}{\Poch{s+k/2}{r}}\,\cR_{k,z}^{(r)}K_s^{k,k},\label{eq:raising-seed}\\
 K_s^{k-2r,k}&=\frac{(-1)^r}{\Poch{s-k/2}{r}}\,\cL_{k,z}^{(r)}K_s^{k,k}.\label{eq:lowering-seed}
\end{align}
The diagonal source coefficient of \(K_s^{t,k}\) is
\begin{equation} c_{t,k}(s)=4\pi\nK(-1)^r\frac{\Gamma(2s)}{\Gamma(s-\ell_-/2)\Gamma(s+\ell_+/2)}. \label{eq:flux-coefficient}\end{equation}
\end{proposition}

\begin{proof}
Fay's mixed kernels are obtained by applying the corresponding Maass operators in either variable; see \cite[Theorem~2.1 and Equations~(45)--(46)]{Fay77}. For integral indices, an analogue in Fay's automorphy-factor convention and variable ordering is recorded in \cite[Theorem~2.6]{BringmannJenkinsKane19}. To fix the normalization for arbitrary real \(k\), differentiate \eqref{eq:explicit-seed} once. The hypergeometric differentiation and contiguous identities \cite[Equations~15.5.3 (with \(n=1\)) and~15.5.20]{DLMF15} give
\begin{equation} \begin{split} \cR_{k,z}K_s^{k,k}&=-\left(s+\frac{k}{2}\right)K_s^{k+2,k},\\ \cL_{k,z}K_s^{k,k}&=-\left(s-\frac{k}{2}\right)K_s^{k-2,k}. \end{split} \label{eq:one-step-mixed-shifts}\end{equation}
All powers used in this differentiation have the branches fixed in \eqref{eq:covariant-factor}, so the calculation is valid for real weights and not only for the integral-weight setting of the second reference. Iteration proves \eqref{eq:raising-seed}--\eqref{eq:lowering-seed}. Multiplying \eqref{eq:fay-identification} by these factors and using \(\Gamma(z+r)=\Poch{z}{r}\Gamma(z)\) gives \eqref{eq:flux-coefficient}. Both sides are meromorphic in \(s\), so the identities extend across the intermediate Pochhammer singularities.
\end{proof}

\begin{remark}[Relation to Fay's construction]
\label{rem:inherited-fay}
For a generic spectral parameter, Fay's Theorem~1.1 and Equation~(13) give a fundamental pair for the radial equation in each angular Fourier mode and supply meromorphic formulas in the parameter \cite[Theorem~1.1 and Equation~(13)]{Fay77}. For the singular \(n=0\) member used as the resolvent seed, the prescribed subdominant asymptotic selects a one-dimensional line, and the logarithmic diagonal coefficient fixes its scalar \cite[Equation~(17) and the proof of Theorem~1.4]{Fay77}. Fay's Theorem~2.1 and Equations~(45)--(46) generate the mixed \(K\)-types from this normalized solution \cite[Theorem~2.1 and Equations~(45)--(46)]{Fay77}. Covariance and the eigenvalue equations alone do not characterize the kernel. The radial ODE, its hypergeometric solution, and the unequal-weight shifts in \eqref{eq:explicit-seed} are prior results. The role of this section is to fix the normalization needed before studying the finite-part actions below.
\end{remark}

\section{Automorphic kernels and orbifold circular finite parts}
\label{sec:direct-kernels}

\subsection{The kernel in the half-plane of absolute convergence}

For each \(\delta\in\overline\Gam\), choose a lift \(\widetilde\delta\in\Gam\). For \(\Re s>1\), define, off the orbit diagonal,
\begin{equation} \Kdir_s^{t,k}(z,w)=\sum_{\delta\in\overline\Gam}\nu(\widetilde\delta)j_k(\widetilde\delta,w)K_s^{t,k}(z,\widetilde\delta w). \label{eq:direct-periodization}\end{equation}
The summand is independent of the choice of \(\widetilde\delta\). The kernel has weight \(t\) in \(z\) and dual weight \(-k\) in \(w\). Away from the orbit diagonal it satisfies
\begin{equation} (\Delta_{t,z}-\lam(s))\Kdir_s^{t,k}=0,\qquad (\Delta_{-k,w}-\lam(s))\Kdir_s^{t,k}=0. \label{eq:kernel-dual-eigen-equations}\end{equation}

Let \(T=\left(\begin{smallmatrix}1&1\\0&1\end{smallmatrix}\right)\) and write
\[\nu(T)=e^{2\pi i\varkappa},\qquad 0\leq\varkappa<1.\]

\begin{theorem}[Kernel in the half-plane of absolute convergence]
\label{thm:direct-kernel}
On compact subsets of
\[\{(s,z,w):\Re s>1,\ z\notin\Gam w\},\]
the series \eqref{eq:direct-periodization} and every fixed coordinate derivative in \(z,w\) converge normally. If \(S\Subset\{\Re s>1\}\) and \(z\) is fixed, then there are constants \(C_{z,S}\) and \(Y_{z,S}\) such that
\begin{equation} \left|\Kdir_s^{t,k}(z,u+iv)\right|\leq C_{z,S}v^{1-\Re s}\qquad(s\in S,\ u\in\Rr,\ v\geq Y_{z,S}). \label{eq:direct-cusp-bound}\end{equation}
For \(\Re s>1\) with \(\lam(s)\) in the resolvent set, let \(\mathscr G_k(s;z,w)\) be the automorphic kernel of \(\Resol_k(s)\). Then
\begin{equation} \Kdir_s^{t,k}(z,w)=c_{t,k}(s)\Dkt{}_z\mathscr G_k(s;z,w) \label{eq:direct-resolvent-kernel}\end{equation}
as an identity of distribution kernels. The product on the right has a unique holomorphic continuation throughout \(\Re s>1\), including the points at which \(\Resol_k(s)\) has a pole, and that continuation equals the left-hand side.
\end{theorem}

\begin{proof}
Put \(d=d_{\Hh}(z,w)\). Since
\[U(z,w)=\sinh^2(d/2),\qquad \sigma(z,w)=\cosh^2(d/2),\]
the covariant factor in \eqref{eq:explicit-seed} satisfies
\begin{equation} \left|P_{t,k}(z,w)U(z,w)^{t/2-s}\right|=\left(\frac{U}{\sigma}\right)^{(t+k)/4}U^{-\Re s}. \label{eq:mixed-seed-modulus}\end{equation}
For \(d\geq1\), the first factor on the right is bounded, including when its exponent is negative. The hypergeometric expansion at \(U^{-1}=0\) and differentiation of \eqref{eq:explicit-seed} give the following estimate without division by a shift coefficient. Let \(S\Subset\{\Re s>1\}\), let \(A_z,B_w\) be fixed invariant differential operators, and choose
\[1<\sigma_0<\sigma_1<\inf_{s\in S}\Re s.\]
For \(d(z,w)\geq1\),
\begin{equation} |A_zB_wK_s^{t,k}(z,w)|\leq C_{S,A,B}(1+d(z,w))^N e^{-\sigma_1d(z,w)} \label{eq:seed-distance-estimate}\end{equation}
for some \(N=N(A,B)\).

For every \(\delta\in\overline\Gam\), one has
\[\frac{4\,\Im z\,\Im(\delta w)}{|\delta w-\overline z|^2}=\cosh^{-2}\!\left(\frac{d(z,\delta w)}{2}\right).\]
Fay's convergence theorem \cite[Theorem~1.5 and Corollary~1.6]{Fay77} and
\[\cosh^{-2}(d/2)\asymp e^{-d}\qquad(d\geq1)\]
imply that
\[\sum_{\delta\in\overline\Gam}e^{-\sigma_0d(z,\delta w)}\]
converges locally uniformly on compact subsets disjoint from the orbit diagonal. Also,
\[(1+d)^Ne^{-\sigma_1d}\leq C_{N,\sigma_0,\sigma_1}e^{-\sigma_0d}\qquad(d\geq1).\]
On such a compact set, the finitely many terms with \(d(z,\delta w)<1\) are smooth. It follows from \eqref{eq:seed-distance-estimate} that the series and all its invariant derivatives converge normally throughout the stated half-plane, including the values at which the Pochhammer factors in \cref{prop:normalized-shifts} vanish. Coordinate and invariant derivatives are comparable on compact charts, so the assertion follows for every fixed coordinate derivative.

Suppose that \(\Re s>1\), that \(\lam(s)\) lies in the resolvent set, that the Pochhammer factors in \cref{prop:normalized-shifts} do not vanish, and that
\[s\notin\left\{\frac{k}{2}-n,-\frac{k}{2}-n:n=0,1,2,\ldots\right\}.\]
For equal weights, \eqref{eq:fay-identification} and \cite[Equation~(16)]{KaraEtAl21} give
\[\Kdir_s^{k,k}(z,w)=c_{k,k}(s)\mathscr G_k(s;z,w)\]
off the diagonal. The sum over \(\overline\Gam\) agrees with the half-sum over \(\Gam\) used there. Applying the shift identities term by term and using \eqref{eq:flux-coefficient} proves \eqref{eq:direct-resolvent-kernel} off the diagonal on this generic set. Both sides are meromorphic in \(s\), while the normally convergent series is holomorphic there. The identity theorem extends the off-diagonal identity across the excluded values of \(s\) at points where \(\lam(s)\) lies in the resolvent set.

At the same generic parameters, Fay's rectangular Fourier expansion \cite[Equation~(71), Theorem~3.1, and Equations~(73) and~(73)']{Fay77} gives a \(v^{1-s}\) zero mode when \(\varkappa=0\), while no zero mode occurs when \(\varkappa\ne0\); in both cases all nonzero modes decrease exponentially for sufficiently large \(v\). Applying the finite-order operator \(\Dkt{}_z\) does not change their \(v\)-exponents. The Fourier coefficients and the normally convergent series are holomorphic in \(s\), so the same expansion extends across the excluded values of \(s\). This proves \eqref{eq:direct-cusp-bound} there as well.

It remains to extend the identity across the orbit diagonal. The equal-weight kernel is normalized by its logarithmic source. By \cref{lem:no-diagonal-counterterm}, a pure raising or lowering string in the first variable gives the circular extension of its off-diagonal expression and introduces no separate diagonal-supported summand. Hence \eqref{eq:direct-resolvent-kernel} holds as an identity of distribution kernels at points \(s\) with \(\Re s>1\) and \(\lam(s)\) in the resolvent set. Normal convergence makes the kernel holomorphic throughout \(\Re s>1\). Near an isolated eigenparameter, the spectral resolvent is meromorphic with a finite-rank principal part; equality on the punctured neighborhood and the identity theorem force every pole of \(c_{t,k}(s)\Dkt\mathscr G_k(s)\) there to be removable, with value equal to the kernel. The same argument extends the identity across the values at which the Pochhammer factors in \cref{prop:normalized-shifts} vanish. This proves the asserted continuation and its uniqueness.
\end{proof}

\subsection{The circular finite part}

For \(r\geq2\), the mixed kernel is generally not locally integrable on the diagonal. For \(r=0,1\) it is locally integrable, but we retain the same circular prescription so that the family has one normalization in all shifts. We use the hyperbolic circular cutoff, including at the two elliptic points of the modular orbifold.

Let \(\pi:\Hh\to X\) be the quotient map.  For a fixed \(z\in X\), choose a lift \(\widetilde z\in\Hh\) and put
\[H_z=\operatorname{Stab}_{\overline\Gam}(\widetilde z),\qquad e_z=|H_z|.\]
For a compatible local lift \(\widetilde w\) of \(w\), use the adapted coordinate
\begin{equation} X_{\widetilde z}(\widetilde w)=\frac{\widetilde w-\widetilde z}{\widetilde w-\overline{\widetilde z}}. \label{eq:adapted-coordinate}\end{equation}
If \(\gamma=\left(\begin{smallmatrix}a&b\\c&d\end{smallmatrix}\right)\in\Gam\), then
\[X_{\gamma\widetilde z}(\gamma\widetilde w)=\frac{c\overline{\widetilde z}+d}{c\widetilde z+d}X_{\widetilde z}(\widetilde w),\]
and the prefactor has modulus one.  Thus the circular cutoff is independent of the compatible choice of lifts.  Within the chosen chart, abbreviate \(X_{\widetilde z}(\widetilde w)\) to \(X_z(w)\).  For sufficiently small \(\rho>0\), define
\[\widetilde q_z:\mathbb D_\rho\longrightarrow X,\qquad \widetilde q_z(X)=\pi\!\left(\frac{\widetilde z-\overline{\widetilde z}X}{1-X}\right).\]
The induced action of \(H_z\) on \(\mathbb D_\rho\) is by rotations, and \(\widetilde q_z\) is \(H_z\)-invariant and induces a homeomorphism from \(\mathbb D_\rho/H_z\) onto its image.  Suppose that an invariant density \(\Omega\), smooth off \(z\), has a polyhomogeneous excised-integral expansion as \(\epsilon\downarrow0\). Write \(\FP_{\epsilon\downarrow0}\) for the constant term of that expansion and define
\begin{equation}
\begin{split}
\FP_z^{\mathrm{circ}}\int_X\Omega:=\;&\frac1{e_z}\FP_{\epsilon\downarrow0}\int_{\epsilon<|X|<\rho}\widetilde q_z^*\Omega\\
&+\int_{X\setminus \widetilde q_z(\mathbb D_\rho)}\Omega.
\end{split}
\label{eq:orbifold-circular-fp}
\end{equation}
We use this definition only for the local types in \cref{lem:circular-fp}; there the indicated constant term is the actual limit for a circular cutoff. Changing \(\rho\) does not change the result. At a regular point this is the full circular finite part; at an elliptic point the factor \(e_z^{-1}\) is the orbifold volume factor.

\begin{lemma}[Fixed-base circular finite part]
\label{lem:circular-fp}
Suppose that, in a smooth local trivialization and after including the smooth measure and test function, a density is a finite sum of terms
\begin{equation}
\begin{aligned}
&X^{-j}A_j(X,\bar X),\quad \bar X^{-j}B_j(X,\bar X) &&(1\leq j\leq J),\\
&C(X,\bar X)\log|X|^2,\quad D(X,\bar X).&&
\end{aligned}
\label{eq:local-fp-types}
\end{equation}
with smooth coefficients. Then the finite part in \eqref{eq:orbifold-circular-fp} exists, and the resulting functional is a distribution of finite order. If the coefficients are meromorphic in a parameter, the pole-cleared finite part is holomorphic, and hence the original family is meromorphic. Here a local pole-clearing factor means a holomorphic function vanishing to at least the order of every pole in the chosen compact parameter disk.
\end{lemma}

\begin{proof}
Fix \(j\). From \(A_j\) subtract its Taylor polynomial of total degree at most \(j-1\). A monomial \(X^a\bar X^b\) in that polynomial would have nonzero angular mean after multiplication by \(X^{-j}\) only if \(a-b=j\), which is impossible when \(a+b\leq j-1\). Thus the polynomial integrates to zero on every circular annulus. The remainder is \(O(|X|^j)\), so the resulting density is locally bounded. The \(\bar X^{-j}\) term is identical, and the logarithmic term is locally integrable. Taylor's remainder estimate bounds the functional by finitely many seminorms of the test function, proving continuity and finite order. Averaging over the rotations in \(H_z\) and dividing by \(e_z\) proves the quotient statement. The same Taylor subtraction is holomorphic coefficientwise after the parameter poles are cleared; removing the pole-clearing factor gives the asserted meromorphic dependence.
\end{proof}

\begin{lemma}[No additional diagonal-supported term]
\label{lem:no-diagonal-counterterm}
Let \(G(z,w)\) be a local Green kernel for the weight Laplacian, with its logarithmic singularity normalized distributionally. Let \(\mathcal D_z\) be either a string of Maass raising operators or a string of Maass lowering operators in the first variable. At fixed \(z\), the distribution \(\mathcal D_zG(z,\cdot)\) equals the \(X_z\)-circular extension of its off-diagonal expression. Thus no separate distribution supported on \(w=z\) is added to that expression. This assertion concerns the kernel distribution; the source equation itself still contains the prescribed delta distribution. The same statement holds in an orbifold chart, with the quotient action normalized by the stabilizer order.
\end{lemma}

\begin{proof}
Write \(X=X_z(w)\) and work first on a uniformizing disk. With \(\partial_X=(\partial_{\Re X}-i\partial_{\Im X})/2\), integration by parts on \(\{|X|>\epsilon\}\), followed by the Taylor subtraction in \cref{lem:circular-fp}, gives
\begin{align}
\partial_X^j\log|X|^2&=(-1)^{j-1}(j-1)!\,\FP^{\mathrm{circ}}\frac1{X^j}, \label{eq:pure-log-derivative}\\
\partial_{\bar X}\!\left(\FP^{\mathrm{circ}}\frac1{X^j}\right)&=\frac{\pi(-1)^{j-1}}{(j-1)!}\partial_X^{j-1}\delta_0,\qquad j\in\Zz_{\geq1}. \label{eq:cauchy-contact-identity}
\end{align}
The conjugate identities hold with \(X\) and \(\bar X\) interchanged. Here \(\FP^{\mathrm{circ}}X^{-j}\) denotes the distribution obtained by integrating over \(\{\epsilon<|X|<\rho\}\) and taking the limit as \(\epsilon\downarrow0\), with the smooth outer contribution understood. The boundary integrals in this calculation are the angular Fourier coefficients removed by the circular Taylor subtraction.

Because the center of \(X_z\) moves with \(z\), derivatives at fixed \(w\) become
\begin{align}
2iy\partial_z&=2iy(\partial_z)_{X,\bar X}-(1-X)\partial_X-\bar X(1-\bar X)\partial_{\bar X}, \label{eq:raising-disc-coordinate}\\
-2iy\partial_{\bar z}&=-2iy(\partial_{\bar z})_{X,\bar X}-X(1-X)\partial_X-(1-\bar X)\partial_{\bar X}. \label{eq:lowering-disc-coordinate}
\end{align}
For \(q\in\Zz_{\geq0}\), let \(\mathscr A_q^+\) be the class of local distributions of the form
\[a_0+a_{\log}\log|X|^2+\sum_{j=1}^q a_j\FP^{\mathrm{circ}}\frac1{X^j},\qquad a_0,a_{\log},a_j\in C^\infty.\]
The first raising operator sends the logarithmic Green expansion into \(\mathscr A_1^+\). Suppose inductively that the result lies in \(\mathscr A_q^+\). The \(\partial_X\)-part of \eqref{eq:raising-disc-coordinate} only increases the order of a pure \(X\)-pole. The only possible contact term comes from the \(\partial_{\bar X}\)-part and \eqref{eq:cauchy-contact-identity}, but its coefficient contains \(\bar X(1-\bar X)\). For every smooth \(a\) and every \(n\in\Zz_{\geq0}\),
\begin{equation} \bar X\,a(X,\bar X)\partial_X^n\delta_0=0. \label{eq:raising-contact-vanishing}\end{equation}
because \(\partial_X\) never differentiates the factor \(\bar X\). This argument also covers all terms produced by the product rule. The base derivative \((\partial_z)_{X,\bar X}\) and the zero-order connection term preserve the same class. This proves the raising induction.

For a lowering string, use the conjugate class
\[\mathscr A_q^-=\left\{a_0+a_{\log}\log|X|^2+\sum_{j=1}^q a_j\FP^{\mathrm{circ}}\frac1{\bar X^j}\right\}.\]
Now every possible contact term is annihilated by \(X\,a(X,\bar X)\partial_{\bar X}^n\delta_0=0\). Hence a pure lowering string also agrees with the circular extension of its off-diagonal expression.

A local Green parametrix has, to arbitrary finite order, a smooth coefficient times \(\log|X|^2\), a smooth term, and a remainder of arbitrarily high regularity. Choosing the expansion order larger than the length of the Maass string makes the preceding induction applicable to the remainder. Finally, the lifted density at an elliptic point is invariant under the finite rotation group \(H_z\). Averaging the identity and dividing its action by \(|H_z|\) gives the normalization in \eqref{eq:orbifold-circular-fp}.
\end{proof}

In a lifted local chart, Fay's angular expansion and singular asymptotics \cite[Theorem~1.1 and Equations~(17)--(19)]{Fay77}, together with the mixed-shift identity \cite[Theorem~2.1 and Equation~(45)]{Fay77}, show that \(\Kdir_s^{t,k}(z,w)f(w)d\mu(w)\) has a local expansion consisting of the terms in \eqref{eq:local-fp-types}. Accordingly, for each fixed base point we make the following definition.

\begin{definition}[Direct finite-part operator]
\label{def:direct-fp}
For \(f\in C_{\mathrm{c},k}^\infty(X,\nu)\), \(\Re s>1\), and fixed \(z\in X\), put
\begin{equation} \cU_{s,\mathrm{dir}}^{t,k}f(z)=\FP_z^{\mathrm{circ}}\int_X\Kdir_s^{t,k}(z,w)f(w)d\mu(w). \label{eq:direct-fp-operator}\end{equation}
For \(\Re s>1\) the resulting distribution depends holomorphically on \(s\); the notation anticipates its meromorphic continuation below.
\end{definition}

\begin{proposition}[Direct realization and source identities]
\label{prop:direct-realization}
For compactly supported smooth inputs and \(\lam(s)\) in the resolvent set,
\begin{equation} \cU_{s,\mathrm{dir}}^{t,k}=c_{t,k}(s)\Dkt\Resol_k(s) \label{eq:direct-operator-identity}\end{equation}
when \(\Re s>1\). The product on the right extends holomorphically across its removable singularities in this half-plane, and the extended identity holds there. Consequently,
\begin{equation}
\begin{split}
(\Delta_t-\lam(s))\cU_{s,\mathrm{dir}}^{t,k}&=c_{t,k}(s)\Dkt,\\
\cU_{s,\mathrm{dir}}^{t,k}(\Delta_k-\lam(s))&=c_{t,k}(s)\Dkt.
\end{split}
\label{eq:source-identities}
\end{equation}
\end{proposition}

\begin{proof}
For \(\Re s>1\) with \(\lam(s)\) in the resolvent set, \cref{thm:direct-kernel} identifies the kernel in \eqref{eq:direct-fp-operator} with the distribution kernel of the right-hand side of \eqref{eq:direct-operator-identity}. By \cref{lem:no-diagonal-counterterm}, the circular prescription in \cref{lem:circular-fp} agrees with its local distributional action. This proves \eqref{eq:direct-operator-identity} there, and holomorphic continuation gives it throughout \(\Re s>1\). The resolvent equation and \(\Delta_t\Dkt=\Dkt\Delta_k\), followed by the same continuation, give \eqref{eq:source-identities}.
\end{proof}

\begin{remark}[Dependence on the cutoff geometry]
\label{rem:shape-dependence}
The cutoff shape is part of the definition. Let \(dA_X\) denote Euclidean area in the disk coordinate. In the local model \(X^{-2}dA_X\), replace \(|X|=\epsilon\) by \(|X|=\epsilon a(\theta)\). The contribution of a constant test function changes by
\[-\int_0^{2\pi}e^{-2i\theta}\log a(\theta)\,d\theta.\]
For \(a(\theta)=e^{\eta\cos2\theta}\), the change is \(-\pi\eta\). Thus replacing the circle by a noncircular cutoff can change the extension by a distribution supported on the diagonal. The hyperbolic circle fixes this ambiguity.
\end{remark}

\begin{remark}[Fixed-base formulation]
\label{rem:no-moving-pv}
The construction does not require differentiation of a moving excised fundamental domain in \(z\). Near an elliptic fixed point, distinct orbit singularities collide and a uniform isolation radius tends to zero. Instead, \eqref{eq:orbifold-circular-fp} is defined at fixed \(z\), while smooth dependence of the output is supplied by the right-hand side of \eqref{eq:direct-operator-identity}. Thus smooth dependence is obtained without invoking a separate moving-cutoff regularity theorem.
\end{remark}

\section{Meromorphic continuation and the cusp coefficient}
\label{sec:continuation}

Weighted spectral and resolvent-kernel constructions go back to Roelcke and Fay; useful accounts of these constructions are \cite{Fay77,KaraEtAl21,Roelcke66,Roelcke67}.

\begin{theorem}[Weighted resolvent continuation]
\label{thm:resolvent-continuation}
The family \(\Resol_k(s)=(\Delta_k-\lam(s))^{-1}\), initially defined for \(\Re s>1/2\) when \(\lam(s)\) lies in the resolvent set, has a meromorphic continuation
\begin{equation} \Resol_k(s):C_{\mathrm{c},k}^\infty(X,\nu)\longrightarrow C_{\mathrm{loc},k}^\infty(X,\nu),\qquad s\in\Cc, \label{eq:resolvent-continuation}\end{equation}
whose principal parts have finite rank. For \(\Re s>1/2\) with \(\lam(s)\) in the resolvent set, the continued family equals the \(L^2\)-resolvent. If the cusp character is nontrivial, \(\Delta_k\) has compact resolvent; if it is trivial, continuation across continuous spectrum is understood in the compact-to-local sense of \eqref{eq:resolvent-continuation}.
\end{theorem}

\begin{proof}
Decompose a cusp section into Fourier modes \(e^{2\pi i(n+\varkappa)x}\). If \(\varkappa\ne0\), their frequencies are bounded away from zero. The resulting cusp estimate makes the embedding of the graph domain of \(\Delta_k\) into \(L_k^2\) compact, so the self-adjoint operator has discrete spectrum and compact resolvent. Its resolvent is then meromorphic with finite-rank principal parts by the spectral theorem. In the real-weight spectral expansion of \cite[Section~2.4, Theorem~5]{KaraEtAl21}, the continuous term is indexed by cusp-fixed vectors; in the present scalar one-cusp setting that term is absent when \(\varkappa\ne0\).

If \(\varkappa=0\), separate the zero Fourier mode. Fix a small parameter neighborhood \(V\) and choose a cusp weight \(\delta\) away from the real parts of the indicial exponents \(s\) and \(1-s\) for \(s\in V\). On the weighted cusp Sobolev spaces \(y^\delta H^m\), variation of parameters with the model solutions \(y^s\) and \(y^{1-s}\) gives a meromorphic inverse for the zero mode. The nonzero modes satisfy the compact cusp estimate used above. Patch these inverses to an interior elliptic parametrix. Between the corresponding weighted spaces this gives
\[(\Delta_k-\lam(s))P_V(s)=\Id+K_V(s),\]
The finitely many polar terms of the model inverse have finite-dimensional ranges. By adjoining these ranges to the auxiliary model space---equivalently, by absorbing their principal parts into a finite-rank meromorphic correction of \(P_V(s)\)---we may take \(K_V(s)\) to be a holomorphic family of compact operators on \(V\). Compactness follows because the cutoff commutators are supported in a fixed cusp collar and the local Sobolev inclusion there is compact. Begin with a neighborhood containing a point \(s_0\) such that \(\Re s_0>1/2\) and \(\lam(s_0)\) lies in the resolvent set. Then \(\Id+K_V(s_0)\) is invertible, so the analytic Fredholm theorem gives a meromorphic inverse on that neighborhood. Continue along a chain of overlapping parameter neighborhoods. On each new overlap the previously constructed resolvent supplies an inverse away from its discrete pole set; hence the nowhere-invertible alternative in the analytic Fredholm theorem is excluded. Uniqueness on overlaps follows from the resolvent identity. The weights are local auxiliary choices: after shrinking an overlap, one may choose a single weight avoiding all indicial exponents there, and the two representatives are compared as maps from compactly supported inputs to local Sobolev spaces. Their equality at a regular point, followed by the identity theorem, makes the compact-to-local continuation independent of that choice. Since every point of the parameter plane can be reached by such a finite chain, this proves \eqref{eq:resolvent-continuation}.

The preceding analytic-Fredholm argument gives the compact-to-local one-cusp continuation used below. Its spectral background goes back to Roelcke's analysis of the weighted automorphic eigenvalue problem \cite{Roelcke66,Roelcke67}. Fredholm theory makes every principal part finite rank. Interior elliptic regularity upgrades the pole-cleared compact-to-local family to the stated smooth mapping property.
\end{proof}

We henceforth define the global meromorphic family by
\begin{equation} \cU_s^{t,k}=c_{t,k}(s)\Dkt\Resol_k(s),\qquad \cK_s^{t,k}(z,w)=c_{t,k}(s)\Dkt{}_z\mathscr G_k(s;z,w). \label{eq:continued-mixed-operator}\end{equation}
By \cref{thm:direct-kernel}, this continues the kernel and operator defined for \(\Re s>1\). The source identities in \eqref{eq:source-identities} consequently continue as meromorphic compact-to-local identities.

\begin{lemma}[Meromorphic local form]
\label{lem:meromorphic-local-form}
Fix a base point \(z\) and a compact parameter disk. After multiplication by a holomorphic factor that clears the finitely many parameter poles, the density \(\cK_s^{t,k}(z,w)d\mu(w)\), near \(w=z\), is a finite sum of the local types in \eqref{eq:local-fp-types}, with coefficients holomorphic in \(s\), plus a remainder of arbitrarily prescribed regularity. Consequently its \(X_z\)-circular action on a smooth test function is meromorphic in \(s\).
\end{lemma}

\begin{proof}
Work on nested uniformizing disks \(D_0\Subset D_1\) about a lift of \(z\). For every integer \(N\), the parameter-dependent local elliptic parametrix construction for \(P(s)=\Delta_k-\lam(s)\) gives a properly supported family \(Q_N(s)\), holomorphic in \(s\), such that the kernels of the errors in \(P(s)Q_N(s)-I\) and \(Q_N(s)P(s)-I\) are \(C^N\) on \(D_0\times D_0\). In the coordinate \(X=X_z(w)\), the kernel of \(Q_N(s)\) has a finite two-dimensional elliptic expansion
\begin{equation} \sum_{a+b\leq M}A_{ab}(s;z,X,\bar X)X^a\bar X^b\log|X|^2+B_N(s;z,X,\bar X), \label{eq:local-parametrix-expansion}\end{equation}
where the coefficients are smooth in the spatial variables and holomorphic in \(s\), and \(B_N\) is \(C^N\) and holomorphic in \(s\). The regularity order \(N\) may be taken arbitrarily large.

The difference between the global resolvent kernel and this local parametrix satisfies an elliptic equation with smooth right-hand side on \(D_0\times D_0\), so it is smooth across the diagonal. It is meromorphic in \(s\). Its principal parts have finite rank; the corresponding resonant and coresonant distributions satisfy elliptic equations and hence are smooth in the interior. On a compact parameter disk \(S\), choose a nonzero holomorphic function \(p_S\) that clears these poles and the scalar poles of \(c_{t,k}\). The pole-cleared smooth remainder is then holomorphic with values in \(C^N(D_0\times D_0)\).

Choose \(N\) larger than the order of \(\Dkt\) plus the desired remainder regularity. Apply \(p_S(s)c_{t,k}(s)\Dkt{}_z\) to \eqref{eq:local-parametrix-expansion}. By \cref{lem:no-diagonal-counterterm}, this produces terms of the forms in \eqref{eq:local-fp-types}, with holomorphic coefficients, plus a remainder of the prescribed regularity. Taylor subtraction in \cref{lem:circular-fp} is continuous on each coefficient space, so its value is holomorphic after pole clearing and meromorphic before it.
\end{proof}

Write
\[\nu(T)=e^{2\pi i\varkappa},\qquad 0\leq\varkappa<1,\qquad T=\begin{pmatrix}1&1\\0&1\end{pmatrix}.\]
When \(\varkappa=0\), let the weight-\(k\) Eisenstein series be normalized by
\begin{equation} E_k(z,s)=y^s+\varphi_k(s)y^{1-s}+E_k^{\mathrm{nc}}(z,s),\qquad \Delta_kE_k(\cdot,s)=\lam(s)E_k(\cdot,s). \label{eq:eisenstein-normalization}\end{equation}

\begin{proposition}[Constant term of the mixed kernel]
\label{prop:kernel-constant-term}
If \(\varkappa=0\), then, for \(v\) sufficiently large,
\begin{equation} \int_{-1/2}^{1/2}\cK_s^{t,k}(z,u+iv)\,du=\beta_s(z)v^{1-s}, \label{eq:kernel-constant-term}\end{equation}
where
\begin{equation} \beta_s(z)=\frac{c_{t,k}(s)}{2s-1}\Dkt E_k(z,s). \label{eq:beta-definition}\end{equation}
Both formulas are meromorphic in \(s\).  When \(\varkappa\ne0\), set \(\beta_s(z)=0\).  After the zero mode is removed, the remaining Fourier modes decrease exponentially in \(v\); if \(\varkappa\ne0\), there is no zero mode. The uniform form of these estimates is given in \cref{lem:uniform-cusp-tails}.
\end{proposition}

\begin{proof}
For \(\Re s>1/2\) at points where \(\lam(s)\) lies in the resolvent set, the zero Fourier coefficient of the resolvent kernel in the second variable solves the Euler equation with solutions \(v^s\) and \(v^{1-s}\). Square integrability removes \(v^s\). In the output/input ordering used here, Fay's symmetry \cite[Equation~(38)]{Fay77} identifies the relevant Fay kernel as the one of index \(-k/2\). Hence Fay's rectangular expansion \cite[Theorem~3.1 and Equation~(73)']{Fay77} has coefficient \(E_{k/2}^{\mathrm{Fay}}(z,s)/(1-2s)\). Under the weight dictionary, \(E_{k/2}^{\mathrm{Fay}}\) is the present \(E_k\), while the resolvent sign in \eqref{eq:fay-laplacian-sign} changes \((1-2s)^{-1}\) into \((2s-1)^{-1}\). Thus the zero mode of \(\mathscr G_k\) is \(E_k(z,s)v^{1-s}/(2s-1)\). Meromorphic continuation proves this identity for every \(s\). Applying \(c_{t,k}(s)\Dkt\) in \(z\) gives \eqref{eq:kernel-constant-term}--\eqref{eq:beta-definition}. The remaining terms are the nonzero Whittaker modes in Equation~(73), which decrease exponentially by Equation~(71) \cite[Equations~(71) and~(73)]{Fay77}.
\end{proof}

\section{Mellin finite-part action on moderate-growth eigenfunctions}
\label{sec:regularized-action}

Let \(\phi\) be a smooth weight-\(k\) automorphic eigenfunction with multiplier \(\nu\),
\begin{equation} \Delta_k\phi=\lambda_\phi\phi,\qquad \phi(x+iy)=O(y^C), \label{eq:moderate-eigenfunction}\end{equation}
for some \(C\), uniformly in a cusp strip.

\begin{lemma}[Constant term of the input]
\label{lem:input-constant-term}
If \(\varkappa\ne0\), set \(\phi^{\mathrm{nc}}=\phi\); then \(\phi^{\mathrm{nc}}\), together with every fixed derivative, decreases exponentially at the cusp. If \(\varkappa=0\), choose \(\sigma_\phi\) with \(\lambda_\phi=\sigma_\phi(1-\sigma_\phi)\). In this case, for \(\sigma_\phi\ne1/2\),
\begin{equation} \phi(x+iy)=A_\phi y^{\sigma_\phi}+B_\phi y^{1-\sigma_\phi}+\phi^{\mathrm{nc}}(x+iy), \label{eq:eigenfunction-constant-term}\end{equation}
where the nonconstant part and its derivatives decrease exponentially. For \(\sigma_\phi=1/2\),
\begin{equation} \phi(x+iy)=y^{1/2}(A_\phi+B_\phi\log y)+\phi^{\mathrm{nc}}(x+iy). \label{eq:eigenfunction-threshold-term}\end{equation}
\end{lemma}

\begin{proof}
The Fourier frequencies are \(n+\varkappa\). Every nonzero frequency satisfies a Whittaker equation; polynomial growth excludes its increasing solution. The zero frequency, present only when \(\varkappa=0\), satisfies \(-y^2f''=\lambda_\phi f\), whose solutions are those displayed above. Their differentiated exponential estimates and the summation over frequencies are established in the next lemma.
\end{proof}

\begin{lemma}[Uniform nonconstant cusp tails]
\label{lem:uniform-cusp-tails}
Put \[\eta_\varkappa=\begin{cases}1,&\varkappa=0,\\ \min\{\varkappa,1-\varkappa\},&0<\varkappa<1.\end{cases}\]
There exists \(c_\varkappa>0\) such that, for every \(a,b\in\Zz_{\geq0}\), there are \(Y_{\phi,a,b},C_{a,b}>0\) such that
\begin{equation}
\left|\partial_x^a(y\partial_y)^b\phi^{\mathrm{nc}}(x+iy)\right|
\leq C_{a,b}e^{-c_\varkappa y}
\qquad(y\geq Y_{\phi,a,b}).
\label{eq:input-uniform-cusp-decay}
\end{equation}
Let \(S\Subset\Cc\), and let \(p_S(s)\) be holomorphic and clear the poles of \(\cK_s^{t,k}\) on a neighborhood of \(S\). For every compact set \(Z\Subset\Hh\), there is a height \(Y_{S,Z}>0\) such that, for all \(q,a,b\in\Zz_{\geq0}\), there is a constant \(C_{S,Z,q,a,b}>0\) satisfying
\begin{equation}
\begin{split}
\bigg|\partial_s^q\partial_u^a(v\partial_v)^b\bigg[p_S(s)\bigg(
&\cK_s^{t,k}(z,u+iv)\\
&-\mathbf 1_{\{\varkappa=0\}}\beta_s(z)v^{1-s}\bigg)\bigg]\bigg|
\leq C_{S,Z,q,a,b}e^{-c_\varkappa v}.
\end{split}
\label{eq:kernel-uniform-cusp-decay}
\end{equation}
This holds for \(z\in Z\), \(s\in S\), and \(v\geq Y_{S,Z}\). The same estimate holds after any fixed derivative in \(z\).
\end{lemma}

\begin{proof}
The multiplier conditions give Fourier expansions
\begin{align*}
\phi(u+iv)&=\sum_{n\in\Zz}a_n(v)e^{2\pi i(n+\varkappa)u},\\
\cK_s^{t,k}(z,u+iv)&=\sum_{m\in\Zz}b_m(s,z;v)e^{2\pi i(m-\varkappa)u}.
\end{align*}
Thus the input and kernel have dual cusp characters, and a zero frequency occurs only for \(\varkappa=0\). Every nonzero frequency \(\xi\) has \(|\xi|\geq\eta_\varkappa\).

The radial equation at frequency \(\xi\ne0\) is a Whittaker equation with argument \(x=4\pi|\xi|v\). A fundamental pair adapted to infinity is described in \cite[Sections~13.14(iv)--(v) and Equation~13.19.3]{DLMF15}, and moderate growth leaves only the decreasing \(W\)-solution. The expansion
\[W_{\alpha,\mu}(x)=e^{-x/2}x^\alpha\bigl(1+O(x^{-1})\bigr)\]
is \cite[Equation~13.19.3]{DLMF15}. The error-bound discussion following that equation, parameter holomorphy \cite[Section~13.14(ii)]{DLMF15}, and compactness of the pole-cleared parameter set make the estimate locally uniform in \((\alpha,\mu)\). Together with the differentiation formulas \cite[Equations~13.15.21--13.15.26]{DLMF15}, this gives, for every \(j\in\Zz_{\geq0}\),
\begin{equation} \partial_v^jW_{\alpha,\mu}(4\pi|\xi|v)\ll (1+|\xi|v)^A e^{-2\pi|\xi|v}, \label{eq:whittaker-uniform-bound}\end{equation}
where chain-rule powers of \(|\xi|\) have been absorbed into the polynomial factor, and \(A\) and the implied constant are uniform on the parameter set.

For the input estimate, fix a sufficiently large height \(Y_0\) and put
\(Y_{\phi,a,b}=2Y_0\). For the kernel estimate, after \(S\) and \(Z\)
are fixed, choose \(Y_0\) above the cusp threshold and so large that
\[Y_0> \sup_{\substack{z\in Z\\ \gamma\in\Gam}} \Im(\gamma z),\]
and put \(Y_{S,Z}=2Y_0\). Thus the horizontal traces used below are smooth and lie above every orbit-diagonal point associated with \(Z\).
 Every decreasing radial coefficient has the form
\[h_\xi(v)=d_\xi W_{\alpha,\mu}(4\pi|\xi|v).\]
When \(|\xi|Y_0\) is sufficiently large, the leading term in Equation~13.19.3 and its differentiated form give
\[|d_\xi|\leq C(1+|\xi|Y_0)^A e^{2\pi|\xi|Y_0}\left(\left|h_\xi(Y_0)\right|+|\xi|^{-1}\left|h_\xi'(Y_0)\right|\right).\]
There are only finitely many remaining frequencies. For each of them, a nonzero Whittaker solution and its derivative cannot vanish simultaneously; ODE uniqueness and compactness of the parameter set therefore give the same inequality after increasing \(C\). Combining this estimate with \eqref{eq:whittaker-uniform-bound} yields
\begin{equation} \left|\partial_v^jh_\xi(v)\right|\leq C(1+|\xi|v)^A e^{-2\pi|\xi|(v-Y_0)}\left(\left|h_\xi(Y_0)\right|+|\xi|^{-1}\left|h_\xi'(Y_0)\right|\right),\quad v\geq2Y_0, \label{eq:whittaker-propagation-bound}\end{equation}
where \(C\) and \(A\) may depend on \(j,Y_0\), and the pole-cleared compact parameter set, but not on \(\xi\) or \(v\).

The untwisted input traces \(e^{-2\pi i\varkappa u}\phi(u+iY_0)\) and \(e^{-2\pi i\varkappa u}\partial_v\phi(u+iv)\big|_{v=Y_0}\) are smooth periodic functions of \(u\).  The untwisted kernel traces \(e^{2\pi i\varkappa u}p_S(s)\cK_s^{t,k}(z,u+iY_0)\) and \(e^{2\pi i\varkappa u}p_S(s)\partial_v\cK_s^{t,k}(z,u+iv)\big|_{v=Y_0}\) are holomorphic in \(s\) with values in \(C^\infty(\Rr/\Zz)\).  Their Fourier coefficients decrease rapidly, uniformly for \(s\in S\) in the kernel case. Summing these propagation estimates over the frequencies proves \eqref{eq:input-uniform-cusp-decay} and \eqref{eq:kernel-uniform-cusp-decay}. Parameter derivatives follow from Cauchy estimates on a slightly larger pole-cleared parameter neighborhood. The same argument applies after fixed derivatives in the base variable.
\end{proof}

For \(Y>1\), let \(X_Y\) be the modular orbifold truncated at height \(Y\). Fix a representative of \(z\) in the standard fundamental domain and take \(Y\) larger than its height. Define the fixed-height action by
\begin{equation} I_{s,Y}(z;\phi)=\FP_z^{\mathrm{circ}}\int_{X_Y}\cK_s^{t,k}(z,w)\phi(w)d\mu(w). \label{eq:truncated-circular-action}\end{equation}
For fixed \(z,Y\), \cref{lem:meromorphic-local-form,lem:circular-fp} make this a meromorphic function of \(s\).

\begin{definition}[Mellin cusp counterterm]
\label{def:cusp-subtraction}
If \(\varkappa\ne0\), put \(Q_{s,\phi}(z;Y)=0\). Suppose \(\varkappa=0\). If \(\sigma_\phi\ne1/2\), set
\begin{equation} Q_{s,\phi}(z;Y)=\beta_s(z)\left(A_\phi\frac{Y^{\sigma_\phi-s}}{\sigma_\phi-s}+B_\phi\frac{Y^{1-\sigma_\phi-s}}{1-\sigma_\phi-s}\right). \label{eq:cusp-subtraction-generic}\end{equation}
If \(\sigma_\phi=1/2\), put \(\delta=1/2-s\) and set
\begin{equation} Q_{s,\phi}(z;Y)=\beta_s(z)Y^\delta\left[\frac{A_\phi}{\delta}+B_\phi\left(\frac{\log Y}{\delta}-\frac1{\delta^2}\right)\right]. \label{eq:cusp-subtraction-threshold}\end{equation}
The right-hand sides are retained as meromorphic functions of \(s\); in particular, \(Y^{\sigma_\phi-s}/(\sigma_\phi-s)\) is not replaced by \(\log Y\) at \(s=\sigma_\phi\). The generic expression is unchanged if \((\sigma_\phi,A_\phi,B_\phi)\) is replaced by \((1-\sigma_\phi,B_\phi,A_\phi)\).
\end{definition}

\begin{definition}[Mellin finite-part transform]
\label{def:regularized-transform}
Away from parameter poles, define
\begin{equation} \cU_{s,\mathrm{fp}}^{t,k}\phi(z)=\lim_{Y\to\infty}\left(I_{s,Y}(z;\phi)-Q_{s,\phi}(z;Y)\right), \label{eq:regularized-action}\end{equation}
provided the limit exists. The construction below gives a meromorphic family. At a genuine pole no value is assigned; only the Laurent expansion is understood.
\end{definition}

The bare symbol \(\cU_s^{t,k}\) denotes the compact-to-local meromorphic operator, while the subscripts \(\mathrm{fp}\) and \(\mathrm{sp}\) denote, respectively, its Mellin finite-part extension to moderate-growth inputs and ordinary integration after specialization.

If a meromorphic vector-valued family \(F_s\) is holomorphic at \(s_0\), or has a removable singularity there, write
\begin{equation} \operatorname{AC}_{s=s_0}F_s:=\widetilde F(s_0), \label{eq:ac-notation}\end{equation}
where \(\widetilde F\) is the holomorphic extension at \(s_0\). This notation assigns no value at a genuine pole. This subtraction follows the general principle of regularized automorphic integration developed by Zagier \cite[Equations~(15), (27), and~(29)--(31)]{Zagier81}, but its coefficient \(\beta_s\) and its parameter dependence are forced here by the mixed resolvent kernel. In a related setting, taking the residue at \(s=\tfrac12\) and letting the cusp-truncation parameter \(T\) tend to infinity need not commute \cite[Section~1.3 and Theorem~2.12(2), (4)]{Wu19}. Here we prescribe only the zero-mode Mellin primitive and compute the resulting defect for the Fay family.

\begin{proposition}[Existence and agreement with the convergent kernel]
\label{prop:regularized-existence}
For every fixed \(z\), \eqref{eq:regularized-action} exists away from a discrete set and extends meromorphically to all \(s\). If \(\Re s>\max\{1,C\}\), then the cusp subtraction tends to zero and
\begin{equation} \cU_{s,\mathrm{fp}}^{t,k}\phi(z)=\FP_z^{\mathrm{circ}}\int_X\Kdir_s^{t,k}(z,w)\phi(w)d\mu(w). \label{eq:direct-agreement}\end{equation}
\end{proposition}

\begin{proof}
Choose a fixed \(Y_0\) above \(z\). Suppose first that \(\varkappa=0\) and \(\sigma_\phi\ne1/2\). Multiplying the constant terms in \cref{prop:kernel-constant-term,lem:input-constant-term}, integrating in the horizontal variable, and including \(dy/y^2\), gives
\[\beta_s(z)\left(A_\phi y^{\sigma_\phi-s-1}+B_\phi y^{-\sigma_\phi-s}\right)dy.\]
Integration gives \eqref{eq:cusp-subtraction-generic}. At the threshold the constant contribution is
\[\beta_s(z)y^{\delta-1}(A_\phi+B_\phi\log y)dy,\]
whose primitive is \eqref{eq:cusp-subtraction-threshold}. If \(\varkappa\ne0\), no constant contribution occurs.

Let \(\phi_0(y)\) denote the zero mode displayed in \eqref{eq:eigenfunction-constant-term} or \eqref{eq:eigenfunction-threshold-term} when \(\varkappa=0\), and put it equal to zero otherwise. Define
\begin{equation} H_s(z,y)=\int_{-1/2}^{1/2}\cK_s^{t,k}(z,u+iy)\phi(u+iy)\,du-\mathbf 1_{\{\varkappa=0\}}\beta_s(z)y^{1-s}\phi_0(y). \label{eq:decaying-horizontal-product}\end{equation}
The Fourier characters of the kernel and input are dual: horizontal integration pairs the frequency \(m-\varkappa\) of the kernel only with the frequency \(-m+\varkappa\) of the input. Therefore \cref{lem:uniform-cusp-tails} gives, after local pole clearing, a common exponentially decreasing majorant for \(H_s\) on every compact parameter set, including after fixed \(s\)- and \(z\)-derivatives. Hence
\begin{equation} \cU_{s,\mathrm{fp}}^{t,k}\phi(z)=I_{s,Y_0}(z;\phi)-Q_{s,\phi}(z;Y_0)+\int_{Y_0}^{\infty}H_s(z,y)\frac{dy}{y^2} \label{eq:fixed-height-formula}\end{equation}
is meromorphic in \(s\). This proves existence without differentiating a moving diagonal cutoff.

If \(\varkappa\ne0\), then \(Q_{s,\phi}(z;Y)=0\), and \eqref{eq:direct-agreement} follows from \cref{thm:direct-kernel,lem:uniform-cusp-tails}.  Assume now that \(\varkappa=0\). If \(A_\phi\ne0\), the growth hypothesis implies \(\Re\sigma_\phi\leq C\), after increasing \(C\) if necessary; if \(B_\phi\ne0\), it similarly implies \(1-\Re\sigma_\phi\leq C\). Thus, in the region \(\Re s>\max\{1,C\}\), every power of \(Y\) occurring in \eqref{eq:cusp-subtraction-generic} tends to zero. The threshold terms also tend to zero because a negative power dominates \(\log Y\). By \cref{thm:direct-kernel}, \(\Kdir_s^{t,k}\) equals the continued kernel in this region, proving \eqref{eq:direct-agreement}.
\end{proof}

\begin{lemma}[Bilinear weighted Green identity]
\label{lem:bilinear-weighted-green}
Let \(D\Subset\Hh\) be a piecewise smooth domain, let \(f\) have weight \(k\), and let \(g\) have dual weight \(-k\). If \(n=(n_x,n_y)\) is the Euclidean outward unit normal and \(ds\) is Euclidean arclength, then
\begin{equation}
\begin{split}
\int_D\left(g\Delta_kf-f\Delta_{-k}g\right)d\mu
=\int_{\partial D}\Big[&\left(f\partial_xg-g\partial_xf+\frac{ik}{y}fg\right)n_x\\
&+\left(f\partial_yg-g\partial_yf\right)n_y\Big]\,ds.
\end{split}
\label{eq:bilinear-weighted-green}
\end{equation}
The contributions on automorphically paired sides cancel. On an upper horizontal boundary \(y=Y\), the contribution is \(\int(f\partial_Yg-g\partial_Yf)\,dx\).
\end{lemma}

\begin{proof}
After multiplication by \(d\mu=dx\,dy/y^2\), one has
\[\frac{g\Delta_kf-f\Delta_{-k}g}{y^2}=\partial_x\left(f\partial_xg-g\partial_xf+\frac{ik}{y}fg\right)+\partial_y\left(f\partial_yg-g\partial_yf\right).\]
The divergence theorem proves \eqref{eq:bilinear-weighted-green}. The automorphy laws show that the displayed Green one-form is invariant; the connection term \(ikfg/y\) is essential here. Opposite boundary orientations therefore make the fluxes on paired sides cancel.
\end{proof}

For the main identity, define the upper-boundary Wronskian
\begin{equation}
\begin{aligned}
W_{s,\phi}(z;Y)=\int_{-1/2}^{1/2}\bigl(&\phi(u+iY)\partial_Y\cK_s^{t,k}(z,u+iY)\\
&-\cK_s^{t,k}(z,u+iY)\partial_Y\phi(u+iY)\bigr)\,du.
\end{aligned}
\label{eq:cusp-wronskian}
\end{equation}

\begin{lemma}[Wronskian cancellation]
\label{lem:wronskian-cancellation}
Locally uniformly after clearing parameter poles,
\begin{equation} W_{s,\phi}(z;Y)-(\lambda_\phi-\lam(s))Q_{s,\phi}(z;Y)\longrightarrow0\qquad(Y\to\infty). \label{eq:cusp-wronskian-cancellation}\end{equation}
\end{lemma}

\begin{proof}
If \(\varkappa\ne0\), all terms decrease exponentially. If \(\varkappa=0\) and \(\sigma_\phi\ne1/2\), insert
\[\cK_{s,0}=\beta_sY^{1-s},\qquad \phi_0=A_\phi Y^{\sigma_\phi}+B_\phi Y^{1-\sigma_\phi}\]
into \eqref{eq:cusp-wronskian}. The result is
\[\beta_s\left[A_\phi(1-\sigma_\phi-s)Y^{\sigma_\phi-s}+B_\phi(\sigma_\phi-s)Y^{1-\sigma_\phi-s}\right].\]
Since
\[\lambda_\phi-\lam(s)=(\sigma_\phi-s)(1-\sigma_\phi-s),\]
this equals \((\lambda_\phi-\lam(s))Q_{s,\phi}\). At \(\sigma_\phi=1/2\), the corresponding expressions are
\[\lambda_\phi-\lam(s)=\delta^2\]
and
\[W_{s,\phi}=\beta_sY^\delta\left[\delta(A_\phi+B_\phi\log Y)-B_\phi\right]=\delta^2Q_{s,\phi}.\]
The nonconstant modes contribute an exponentially small remainder.
\end{proof}

\begin{theorem}[Finite-part spectral reduction]
\label{thm:regularized-reduction}
For every moderate-growth eigenfunction satisfying \eqref{eq:moderate-eigenfunction},
\begin{equation} (\lambda_\phi-\lam(s))\cU_{s,\mathrm{fp}}^{t,k}\phi=c_{t,k}(s)\Dkt\phi. \label{eq:regularized-reduction}\end{equation}
Equivalently, as an identity of meromorphic \(C_{\mathrm{loc},t}^\infty\)-valued families,
\begin{equation} {\cU_{s,\mathrm{fp}}^{t,k}\phi=\frac{c_{t,k}(s)}{\lambda_\phi-s(1-s)}\Dkt\phi.} \label{eq:regularized-scalar-formula}\end{equation}
\end{theorem}

\begin{proof}
Put \(g(w)=\cK_s^{t,k}(z,w)\). In the second variable, \(g\) has dual weight \(-k\) and satisfies, distributionally,
\begin{equation} (\Delta_{-k,w}-\lam(s))g(w)=c_{t,k}(s)\Dkt{}_z\delta_z(w), \label{eq:kernel-second-variable-source}\end{equation}
where the right-hand side is characterized by its action \(f\mapsto c_{t,k}(s)\Dkt f(z)\). This source identity holds first for \(\Re s>1\) when \(\lam(s)\) lies in the resolvent set, by \cref{prop:direct-realization}; both sides then agree meromorphically by \cref{thm:resolvent-continuation}. Apply the distributional form of \cref{lem:bilinear-weighted-green} on \(X_Y\), obtained by testing against smooth cutoffs that approach the characteristic function of \(X_Y\). The paired sides cancel, and the top boundary is \(W_{s,\phi}(z;Y)\). The diagonal pairing with \(g\) is the circular finite part by \cref{lem:no-diagonal-counterterm}. Consequently,
\[\left\langle(\Delta_{-k,w}-\lam(s))g,\phi\right\rangle_{X_Y}=(\lambda_\phi-\lam(s))I_{s,Y}(z;\phi)-W_{s,\phi}(z;Y).\]
The left-hand side is the source action in \eqref{eq:kernel-second-variable-source}. Therefore
\begin{equation} (\lambda_\phi-\lam(s))I_{s,Y}(z;\phi)=c_{t,k}(s)\Dkt\phi(z)+W_{s,\phi}(z;Y). \label{eq:truncated-green}\end{equation}
Subtract \((\lambda_\phi-\lam(s))Q_{s,\phi}\), let \(Y\to\infty\), and use \cref{lem:wronskian-cancellation}. This proves \eqref{eq:regularized-reduction} away from a discrete set and hence meromorphically everywhere. Division wherever \(\lambda_\phi\ne\lam(s)\), followed by meromorphic continuation, proves \eqref{eq:regularized-scalar-formula}. The right-hand side also supplies the claimed smooth dependence on \(z\), so no moving-cutoff regularity theorem is being used.
\end{proof}

\begin{remark}[Normalization conditions]
\label{rem:conditional-canonicity}
The value in \eqref{eq:regularized-action} is fixed by three requirements: it agrees with the integral in \eqref{eq:direct-agreement}; it subtracts a Mellin antiderivative of the kernel's constant Fourier mode; and that antiderivative tends to zero for \(\Re s>\max\{1,C\}\). Adding a \(Y\)-independent meromorphic counterterm would violate the first requirement. This is a conditional normalization statement, not an assertion that every regularization scheme must give the same value.
\end{remark}

The next theorem isolates the mechanism used in the holomorphic example. It shows that the order of specialization and Mellin finite-part continuation can matter even when the specialized integral itself converges.

\begin{theorem}[Boundary defect at a simple source-coefficient zero]
\label{thm:boundary-defect}
Assume \(\varkappa=0\), let \(s_0=\sigma_\phi\ne1/2\), and suppose that \(c_{t,k}\) has a simple zero at \(s_0\) and that \(E_k(z,s)\) is regular there. Assume also that \(\cK_s^{t,k}\) has a removable singularity at \(s_0\) as a distribution kernel, and denote its continued value by \(\cK_{s_0}^{t,k}\). Let \(\phi\) have the constant term
\[A_\phi y^{s_0}+B_\phi y^{1-s_0}.\]
If the integral of the specialized kernel \(\cK_{s_0}^{t,k}\) against \(\phi\) converges, denote it by \(\cU_{s_0,\mathrm{sp}}^{t,k}\phi\). Then
\begin{equation} \operatorname{AC}_{s=s_0}\cU_{s,\mathrm{fp}}^{t,k}\phi-\cU_{s_0,\mathrm{sp}}^{t,k}\phi=A_\phi\,\beta_{s_0}', \label{eq:boundary-defect}\end{equation}
where
\begin{equation} \beta_{s_0}'=\frac{c_{t,k}'(s_0)}{2s_0-1}\Dkt E_k(\cdot,s_0). \label{eq:beta-derivative}\end{equation}
\end{theorem}

\begin{proof}
Because \(\beta_{s_0}=0\), the second term in \eqref{eq:cusp-subtraction-generic} vanishes at \(s_0\), while
\[\operatorname{AC}_{s=s_0}\!\left[\beta_s A_\phi\frac{Y^{s_0-s}}{s_0-s}\right]=-A_\phi\beta_{s_0}'.\]
Thus \(\operatorname{AC}_{s=s_0}Q_{s,\phi}=-A_\phi\beta_{s_0}'\). For a fixed \(Y_0\), \cref{lem:meromorphic-local-form} and the removability hypothesis identify the continued value of \(I_{s,Y_0}\) with the circular action of \(\cK_{s_0}^{t,k}\). In the fixed-height formula \eqref{eq:fixed-height-formula}, \cref{lem:uniform-cusp-tails} permits specialization under the remaining tail integral. Since \(\beta_{s_0}=0\), this is the tail of the specialized-kernel integral. Taking the removable value therefore gives \(\operatorname{AC}_{s=s_0}\cU_{s,\mathrm{fp}}^{t,k}\phi=\cU_{s_0,\mathrm{sp}}^{t,k}\phi-\operatorname{AC}_{s=s_0}Q_{s,\phi}\), proving \eqref{eq:boundary-defect}. Differentiating \eqref{eq:beta-definition} at the simple zero gives \eqref{eq:beta-derivative}.
\end{proof}

\section{Source-coefficient zeros, spectral projections, and cusp defects}
\label{sec:spectral-selectors}

\subsection{An isolated-eigenvalue expansion}

The next result identifies the effect of matching a source-coefficient zero with an isolated spectral value. Its value is a scalar multiple of the Maass shift of the spectral projection, and its first derivative contains the reduced resolvent.

\begin{theorem}[Spectral projection and reduced resolvent]
\label{thm:selector}
Let \(\Re s_0>1/2\) and \(\lambda_0=\lam(s_0)\). Assume that \(\lambda_0\) is either in the resolvent set of \(\Delta_k\) or is an isolated point of its spectrum. In the second case let \(P_0\) be its orthogonal spectral projection; in the first case put \(P_0=0\). Put \(Q_0=\Id-P_0\) and
\begin{equation} S_0=Q_0(\Delta_k-\lambda_0)^{-1}Q_0. \label{eq:reduced-resolvent}\end{equation}
At an isolated eigenvalue the inverse in this formula is taken on \(\operatorname{Ran}Q_0\). Suppose that \(c_{t,k}\) has a simple zero at \(s_0\). Write
\begin{equation} h=s-s_0,\qquad a=2s_0-1,\qquad c_{t,k}(s)=c_1h+c_2h^2+O(h^3), \label{eq:c-expansion}\end{equation}
so \(c_1=c_{t,k}'(s_0)\) and \(c_2=c_{t,k}''(s_0)/2\). As maps from compactly supported smooth vectors to local smooth vectors, \(\cU_s^{t,k}\) is holomorphic at \(s_0\) and
\begin{equation} {\cU_{s_0}^{t,k}=\frac{c_1}{a}\Dkt P_0.} \label{eq:spectral-selector}\end{equation}
Its first derivative is
\begin{equation} {(\cU^{t,k})'_{s_0}=\left(\frac{c_2}{a}-\frac{c_1}{a^2}\right)\Dkt P_0+c_1\Dkt S_0.} \label{eq:selector-derivative}\end{equation}
Consequently,
\begin{equation} {\Dkt S_0=\frac1{c_1}\left[(\cU^{t,k})'_{s_0}-\left(\frac{c_2}{c_1}-\frac1a\right)\cU_{s_0}^{t,k}\right].} \label{eq:reduced-resolvent-recovery}\end{equation}
\end{theorem}

\begin{proof}
For \(\Re s>1/2\) at points where \(\lam(s)\) lies in the resolvent set, the continued family agrees with the \(L^2\)-resolvent. The self-adjoint reduced-resolvent expansion therefore applies near \(s_0\). Substituting \(\lam(s_0+h)-\lambda_0=-h(a+h)\) into the reduced-resolvent expansion \cite[Chapter~III, Section~6.5, especially Equation~(6.32)]{Kato95} gives
\begin{equation} \Resol_k(s_0+h)=\frac{P_0}{h(a+h)}+S_0-ahS_0^2+O(h^2), \label{eq:resolvent-laurent}\end{equation}
because \(\lam(s_0+h)=\lambda_0-ah-h^2\). In the resolvent-set case the term containing \(P_0=0\) is absent and \(S_0\) is the ordinary resolvent at \(\lambda_0\). Multiplying \eqref{eq:resolvent-laurent} by \(c_{t,k}(s_0+h)\Dkt\) gives
\[\cU_{s_0+h}^{t,k}=\frac{c_1}{a}\Dkt P_0+h\left[\left(\frac{c_2}{a}-\frac{c_1}{a^2}\right)\Dkt P_0+c_1\Dkt S_0\right]+O(h^2).\]
The three assertions follow by comparing coefficients.
\end{proof}

\begin{remark}[Thresholds and resonances]
\label{rem:selector-scope}
The assumptions are essential. In the weight-zero cofinite case, at a small isolated eigenvalue \(\lambda_0=s_0(1-s_0)<1/4\), Risager proves the cusp mapping \(R_0(s_0):B_0\to B_{1-s_0+\varepsilon}\) for every \(\varepsilon>0\) \cite[Section~2.4, Proposition~2.3]{Risager11}. At \(s_0=1/2\), the spectral denominator is quadratic in \(h\), so a simple zero of \(c\) does not remove the pole. At an embedded eigenvalue or a resonance, \(Q_0\Resol_k(s)Q_0\) can have an additional principal part and \(S_0\) need not be a bounded reduced inverse. When \(\varkappa=0\) and \(\Re s_0<1/2\), the continued resolvent may likewise have a resonance even when \(\lambda(s_0)\) is in the self-adjoint resolvent set. The theorem makes no claim in these situations. Its use here is the explicit Fay-kernel realization of \(\Dkt P_0\) and \(\Dkt S_0\).
\end{remark}

\subsection{The holomorphic discrete-series point}

Assume for the rest of the paper that the multiplier is trivial, that \(k>2\) is even, and that
\begin{equation} t=k+2r,\qquad s_0=\frac k2,\qquad \nK=\frac{k-1}{4\pi}. \label{eq:holomorphic-normalization}\end{equation}
Let \(P_k^{\mathrm{hol}}\) be the orthogonal projection in the unitary weight-\(k\) model onto
\[\mathcal H_k^{\mathrm{hol}}=\{y^{k/2}f(z):f\in S_k(\Gam)\}.\]
This is the eigenspace with
\begin{equation} \lambda_0=\frac k2\left(1-\frac k2\right). \label{eq:holomorphic-eigenvalue}\end{equation}
In the present trivial-character case, the continuous term in the spectral expansion \cite[Section~2.4, Equation~(10), and Theorem~5]{KaraEtAl21} is parametrized by \(s=1/2+i\tau\). Equation~(10) assigns to it the eigenvalue \(1/4+\tau^2\), so the essential spectrum begins at \(1/4\). Since \(k(2-k)/4<1/4\), this eigenvalue is isolated when the space is nonzero and belongs to the resolvent set when the space is zero. Since \(\cR_{k-2}\) is the negative adjoint of \(\cL_k\), the factorization
\[\cR_{k-2}\cL_k=-\Delta_k-\frac{k(k-2)}4\]
shows that an eigenfunction with eigenvalue \(k(2-k)/4\) lies in \(\ker\cL_k\). This is equivalent to the holomorphy of \(y^{-k/2}\phi\), and square integrability then forces cuspidality.

\begin{corollary}[Maass shift of the Bergman projection]
\label{cor:holomorphic-projector}
Under \eqref{eq:holomorphic-normalization},
\begin{equation} {\cU_{k/2}^{k+2r,k}=\frac{(-1)^r}{\Poch{k}{r}}\cR_k^{(r)}P_k^{\mathrm{hol}}.} \label{eq:holomorphic-projector}\end{equation}
Moreover,
\begin{equation} {\frac1{r!}\cL_{k+2r}^{(r)}\cU_{k/2}^{k+2r,k}=P_k^{\mathrm{hol}}.} \label{eq:lowering-recovers-projector}\end{equation}
For \(r=0\), \eqref{eq:holomorphic-projector} is the ordinary Bergman projection; compare the normalized reproducing-kernel identity in \cite[Definition~7 and Theorem~4(3)]{Steiner16}.
\end{corollary}

\begin{proof}
At \(s_0=k/2\), \eqref{eq:flux-coefficient} becomes
\[c_{k+2r,k}(s)=(k-1)(-1)^r\frac{\Gamma(2s)}{\Gamma(s-k/2)\Gamma(s+k/2+r)}.\]
Since \(1/\Gamma(h)=h+O(h^2)\),
\begin{equation} c_{k+2r,k}'(k/2)=(k-1)(-1)^r\frac{\Gamma(k)}{\Gamma(k+r)}. \label{eq:holomorphic-c-prime}\end{equation}
Now \(2s_0-1=k-1\), so \cref{thm:selector} gives \eqref{eq:holomorphic-projector}. Repeated factorization gives, for every \(\psi\in\mathcal H_k^{\mathrm{hol}}\),
\begin{equation} \cL_{k+2r}^{(r)}\cR_k^{(r)}\psi=(-1)^r r!\Poch{k}{r}\,\psi. \label{eq:reverse-ladder}\end{equation}
Combining this with \eqref{eq:holomorphic-projector} proves \eqref{eq:lowering-recovers-projector}.
\end{proof}

\begin{remark}[Fay's reproducing kernel]
\label{rem:projector-prior-art}
The residue and reproducing-kernel mechanism behind \eqref{eq:holomorphic-projector} is already present in \cite[Equation~(39); Theorem~2.1 and Equations~(45)--(49); Corollary~3.2]{Fay77}. \Cref{cor:holomorphic-projector} fixes its normalization in the unitary weight convention used here; it is not a new Bergman-kernel construction.
\end{remark}

The first derivative also contains the reduced resolvent. Let \(S_k^{\mathrm{hol}}\) be the reduced resolvent at \eqref{eq:holomorphic-eigenvalue} and put
\[H_n=\sum_{j=1}^n\frac1j,\qquad \eta_{k,r}=2H_{k-1}-H_{k+r-1}-\frac1{k-1}.\]
Expanding the gamma quotient one order further gives \(c_2/c_1=2H_{k-1}-H_{k+r-1}\). Hence \cref{thm:selector} yields the reduced-resolvent identity
\begin{equation} {\cR_k^{(r)}S_k^{\mathrm{hol}}=\frac{(-1)^r\Poch{k}{r}}{k-1}\left[(\cU^{k+2r,k})'_{k/2}-\eta_{k,r}\cU_{k/2}^{k+2r,k}\right].} \label{eq:holomorphic-reduced-resolvent}\end{equation}
Thus the value and spectral derivative of the explicit mixed-kernel family together recover the shifted solution operator on the orthogonal complement of the holomorphic cusp space. Spectral derivatives also occur in the construction of polyharmonic Maass forms; see, for comparison, \cite[Theorem~1.1(3); Section~6]{AndersenLagariasRhoades19}.

\subsection{Holomorphic forms and the cusp defect}

Let \(M_k(\Gam)\) and \(S_k(\Gam)\) denote the holomorphic modular and cusp forms; we use the level-one dimension and decomposition results in \cite[Theorem~3.5.2]{DiamondShurman05}. Normalize the holomorphic Eisenstein series by \(E_k^{\mathrm{hol}}(z)=1+O(e^{2\pi iz})\). The unitary Eisenstein series in \eqref{eq:eisenstein-normalization} satisfies
\begin{equation} E_k(z,k/2)=y^{k/2}E_k^{\mathrm{hol}}(z). \label{eq:holomorphic-eisenstein-value}\end{equation}
This follows by evaluating its defining Poincar\'e series at \(s=k/2\).

For \(F\in M_k(\Gam)\), let \(a_0(F)\) be its constant Fourier coefficient and put
\begin{equation} F^\circ=F-a_0(F)E_k^{\mathrm{hol}}\in S_k(\Gam),\qquad \Phi_F(z)=y^{k/2}F(z). \label{eq:holomorphic-decomposition}\end{equation}
Although \(\Phi_F\) need not lie in \(L^2\), its Petersson pairing with every cusp form converges. These pairings extend the finite-rank cusp projection to the present inputs; write this extension as \(P_{k,\mathrm{ext}}^{\mathrm{hol}}\). Unfolding the Eisenstein term, or ordinary cusp--Eisenstein orthogonality, gives
\begin{equation} P_{k,\mathrm{ext}}^{\mathrm{hol}}\Phi_F=\Phi_{F^\circ}. \label{eq:extended-cusp-projection}\end{equation}

Recall that the notation in \eqref{eq:ac-notation} means continuation in the parameter \(s\). In general it is not a regularized integral of the already specialized integrand \(\cK_{s_0}^{k+2r,k}(z,\cdot)\Phi_F\).

\begin{theorem}[Specialization defect for holomorphic forms]
\label{thm:holomorphic-defect}
For \(F\in M_k(\Gam)\), the integral of the specialized kernel is absolutely convergent and
\begin{equation} {\cU_{k/2,\mathrm{sp}}^{k+2r,k}\Phi_F=\frac{(-1)^r}{\Poch{k}{r}}\cR_k^{(r)}\Phi_{F^\circ}.} \label{eq:holomorphic-specialized}\end{equation}
In contrast, analytic continuation of the Mellin finite-part family gives
\begin{equation} {\operatorname{AC}_{s=k/2}\cU_{s,\mathrm{fp}}^{k+2r,k}\Phi_F=\frac{(-1)^r}{\Poch{k}{r}}\cR_k^{(r)}\Phi_F.} \label{eq:holomorphic-ac}\end{equation}
Therefore specialization and analytic continuation of the finite-part family fail to commute by the explicit amount
\begin{equation} {\operatorname{AC}_{s=k/2}\cU_{s,\mathrm{fp}}^{k+2r,k}\Phi_F-\cU_{k/2,\mathrm{sp}}^{k+2r,k}\Phi_F=\frac{(-1)^ra_0(F)}{\Poch{k}{r}}\cR_k^{(r)}\!\left(y^{k/2}E_k^{\mathrm{hol}}\right).} \label{eq:holomorphic-defect}\end{equation}
\end{theorem}

\begin{proof}
At \(s=k/2\), the zero of \(c_{k+2r,k}\) cancels the resolvent pole when the cusp space is nonzero; if that space vanishes, it instead multiplies a regular resolvent. In both cases the resulting kernel is the finite-rank Maass shift of the cusp projection (possibly zero) in \eqref{eq:holomorphic-projector}; it is smooth on the diagonal and exponentially decreasing in the integration variable. Its integral against \(\Phi_F\) therefore converges absolutely and equals the raised extension in \eqref{eq:extended-cusp-projection}, which is \eqref{eq:holomorphic-specialized}.

For the parameter-dependent action, apply \eqref{eq:regularized-scalar-formula} with \(\lambda_\phi=\lam(k/2)\), then take the removable value. Equations \eqref{eq:holomorphic-c-prime} and \(2s_0-1=k-1\) give \eqref{eq:holomorphic-ac}. Subtracting \eqref{eq:holomorphic-specialized} and using \eqref{eq:holomorphic-decomposition} proves \eqref{eq:holomorphic-defect}. Equivalently, this is \cref{thm:boundary-defect} with \(A_{\Phi_F}=a_0(F)\).
\end{proof}

\begin{remark}[Relation to Zagier regularization]
\label{rem:zagier-distinction}
For a fixed integrable automorphic function in the regularizable class, Zagier's renormalized integral agrees with the ordinary integral \cite[Equation~(31) and the sentence immediately following it]{Zagier81}; compare also \cite[Theorem~2.12(4) and Definition~2.13]{Wu19}. Accordingly, applying that fixed-function regularization to \(\cK_{k/2}\Phi_F\) gives the specialized value \eqref{eq:holomorphic-specialized}. The other side of \eqref{eq:holomorphic-defect} is instead the analytic-continuation value of an \(s\)-dependent kernel family. The theorem is a noncommutation statement, not a competing value for a fixed convergent integral.
\end{remark}

\subsection{The weight-twelve calculation}

Let
\[\Phi_\Delta=y^6\Delta_{\mathrm{Ram}},\qquad \Phi_E=y^6E_{12}^{\mathrm{hol}}.\]
Since \(S_{12}(\Gam)=\Cc\Delta_{\mathrm{Ram}}\), write \(P_\Delta\) for the orthogonal projection onto \(\Cc\Phi_\Delta\).

\begin{corollary}[Weight-twelve specialization defect]
\label{cor:weight-twelve-defect}
For \(k=12,t=14,r=1\), one has
\begin{equation} c_{14,12}(s)=-11\frac{\Gamma(2s)}{\Gamma(s-6)\Gamma(s+7)},\qquad c_{14,12}'(6)=-\frac{11}{12}. \label{eq:weight12-c}\end{equation}
Moreover,
\begin{equation} {\cU_6^{14,12}=-\frac1{12}\cR_{12}P_\Delta,\qquad \cL_{14}\cU_6^{14,12}=P_\Delta.} \label{eq:delta-projector}\end{equation}
The specialized integral against \(\Phi_E\) converges, but
\begin{equation} {\cU_{6,\mathrm{sp}}^{14,12}\Phi_E=0.} \label{eq:e12-specialized}\end{equation}
On the other hand,
\begin{equation} {\operatorname{AC}_{s=6}\cU_{s,\mathrm{fp}}^{14,12}\Phi_E=-\frac1{12}\cR_{12}\Phi_E\ne0.} \label{eq:e12-regularized}\end{equation}
\end{corollary}

\begin{proof}
Equation \eqref{eq:weight12-c} follows from \eqref{eq:flux-coefficient}. Since the bottom eigenspace is \(\Cc\Phi_\Delta\), \cref{thm:selector} gives the first identity in \eqref{eq:delta-projector}; the second follows from \(\cL_{14}\cR_{12}=-\Delta_{12}-42\) on the \(-30\)-eigenspace. The specialized kernel is smooth and rapidly decreasing in its second variable, so its ordinary integral realizes \(-\frac1{12}\cR_{12}P_{12,\mathrm{ext}}^{\mathrm{hol}}\). Equation \eqref{eq:extended-cusp-projection} gives \(P_{12,\mathrm{ext}}^{\mathrm{hol}}\Phi_E=0\), which proves \eqref{eq:e12-specialized}. \Cref{thm:holomorphic-defect}, with \(a_0(E_{12}^{\mathrm{hol}})=1\), gives \eqref{eq:e12-regularized}. If \(\cR_{12}\Phi_E\) vanished, then applying \(\cL_{14}\) would contradict \(\cL_{14}\cR_{12}\Phi_E=-12\Phi_E\ne0\); hence the continued value is nonzero.
\end{proof}

The sign is also visible in the cusp expansion. With \(h=s-6\),
\begin{equation} \beta_s=-\frac h{12}\cR_{12}\Phi_E+O(h^2),\qquad Q_{s,\Phi_E}(Y)=\beta_s\frac{Y^{-h}}{-h}. \label{eq:e12-counterterm}\end{equation}
Hence \(\operatorname{AC}_{s=6}Q_{s,\Phi_E}=\cR_{12}\Phi_E/12\). Consequently, the prescription \(I-Q\) has the value in \eqref{eq:e12-regularized}, although the specialized integral is zero.

\begin{remark}[Escaping-mass model]
\label{rem:escape-toy}
The loss of uniform integrability is already visible in
\[f_h(y)=h y^{-1-h},\qquad \int_1^\infty f_h(y)\,dy=1\quad(\Re h>0),\]
whereas \(f_0\equiv0\). For real \(h>0\), under \(u=h\log y\), the measure \(f_h(y)dy\) becomes \(e^{-u}du\): as \(h\to0\), its mass moves to arbitrarily large cusp height. Equation \eqref{eq:e12-counterterm} is the automorphic, mixed-weight version of this elementary loss of uniform integrability.
\end{remark}

\section{Conclusion}
\label{sec:conclusion}

The local hypergeometric and mixed-weight kernels used here are Fay's, and the abstract reduced-resolvent expansion is Kato's. The additional content is the joint realization of the circular diagonal extension and the Mellin cusp subtraction, together with the spectral-parameter derivative formula \eqref{eq:boundary-defect}. The weight-twelve example shows concretely that specializing the kernel before integration and evaluating the parameter-dependent finite-part family after continuation need not agree. Their discrepancy is the explicit cusp coefficient above, not an alternative regularized value of the fixed convergent integral. Extending the defect formula from the modular surface to several cusps would replace the scalar \(\beta_s\) by a scattering-matrix boundary jet; that extension is not claimed here.

\newpage
\section*{Declarations}
{Funding.}
This research received no external funding.

{Competing Interests.}
The author declares no competing interests.

{Use of AI-assisted technologies.}
ChatGPT, Gemini, and Fable 5 were used to assist with proof development and language editing.

\bibliographystyle{amsplain}
\bibliography{maass-fay-reconstructed}

\end{document}